\documentclass[a4paper]{article}
\usepackage{amsmath}
\usepackage{graphics}
\usepackage{latexsym}
\usepackage{amsfonts}
\usepackage{mathrsfs}
\usepackage{amssymb}
\usepackage{amssymb,amsmath,amsthm, amsfonts}
\numberwithin{equation}{section}

\newtheorem{propo}{Proposition}[section]

\newtheorem{corollary}{Corollary}

\newtheorem{remark}{Remark}

\newtheorem{assumption}{Assumption}

\def\qed{ \ \vrule width.2cm height.2cm depth0cm\smallskip}

\def \ind{1\!\!1}
\def \x{X^{t,x}}
\def \ed {\end{document}}
\def \tx {(t,x)\in \esp}

\def \pr {{\bf Proof}}

\def \ij {i\in \cJ}

\def \lb{\label}

\newcommand{\eps}{\varepsilon}

\newcommand{\brm}{\begin{rem}}
\newcommand{\ermq}{\end{rem}}

\newcommand{\ba}{\begin{array}}
\newcommand{\ea}{\end{array}}
\newcommand{\be}{\begin{equation}}
\newcommand{\ee}{\end{equation}}
\newcommand{\bea}{\begin{eqnarray}}
\newcommand{\eea}{\end{eqnarray}}
\newcommand{\beaa}{\begin{eqnarray*}}
\newcommand{\eeaa}{\end{eqnarray*}}

\def \R{I\!\!R}
\def \E{\mathbb{E}}

\def \ij{(i,j)\in \gam}

\def\g{\gamma}
\def\d{\delta}

\def\l{\lambda}

\def\si{\sigma}
\def\t{\tau}

\def\L{\Lambda}

\def\cA{{\cal A}}

\def\cC{{\cal C}}

\def\cF{{\cal F}}

\def\cH{{\cal H}}

\def\cJ{{\cal J}}

\def\cO{{\cal O}}

\def\cS{{\cal S}}

\def\cU{{\cal U}}

\def\no{\noindent}

\def\ms{\medskip}
\def\bs{\bigskip}
\def\q{\quad}
\def\qq{\qquad}

\def\bF{{\bf F}}

\def\qed{ \hfill \vrule width.25cm height.25cm depth0cm\smallskip}

\newcommand{\basa}{\begin{assumption}}
\newcommand{\easa}{\end{assumption}}

\newcommand{\bas}{\begin{assum}}
\newcommand{\eas}{\end{assum}}

\def\limsup{\mathop{\overline{\rm lim}}}
\def\liminf{\mathop{\underline{\rm lim}}}

\def\esssup{\mathop{\rm esssup}}
\def\essinf{\mathop{\rm essinf}}

\def\dis{\displaystyle}

\def\bF{{\bf F}}

\def \P{\mathbb{P}}
\newtheorem{thm}{Theorem}[section]
\newtheorem{lem}[thm]{Lemma}
\newtheorem{prop}[thm]{Proposition}
\newtheorem{rem}[thm]{Remark}

\newtheorem{assum}[thm]{Assumption}

\newtheorem{axiom}{Definition}
\newcommand{\rw}{\rightarrow}
\def \R{\mathbb{R}}

\def \esssup {\mbox{ess sup}}
\def \essinf {\mbox{ess inf}}
\def \esp {[0,T]\times \R^k}
\def \gam {\Gamma^{1}\times \Gamma^{2} }
\def\gami{{(\Gamma^1)}^{-i}}
\def\gamj{{(\Gamma^2)}^{-j}}

\title{Viscosity Solutions of Systems of Variational Inequalities with Interconnected Bilateral Obstacles.}
\author{Boualem Djehiche\footnote{This work was completed while the first author was visiting the Department of Mathematics of Universit\'e du Maine. Financial support from GEANPYL  and Svensk Export Kredit (SEK) is gratefully acknowledged.}\,\,\thanks{Department of Mathematics, Royal Institute of Technology, SE-100 44 Stockholm, Sweden,  boualem@math.kth.se.} ,\, \, Said Hamad\`ene\thanks{Universit\'e du
Maine, LMM, Avenue Olivier Messiaen, 72085 Le Mans, Cedex 9, France, e-mail: hamadene@univ-lemans.fr} \,\, and \, Marie-Amelie Morlais\thanks{Universit\'e du Maine, LMM, Avenue Olivier Messiaen, 72085 Le Mans, Cedex 9, France, e-mail: Marie$_-$Amelie.Morlais@univ-lemans.fr} }

\begin{document}
\date{\today}
\maketitle

\begin{abstract}
We study a general class of nonlinear second-order variational inequalities with interconnected bilateral obstacles, related to a multiple modes switching game. Under rather weak assumptions, using systems  of penalized unilateral backward SDEs, we construct a continuous viscosity solution of polynomial growth. Moreover, we  establish a comparison result which in turn yields uniqueness of the solution.

\end{abstract}
\no{\bf AMS Classification subjects}: 49N70, 49L25, 60H30, 90C39, 93E20
\medskip

\no {$\bf Keywords$}: Switching games, variational inequalities,
oblique reflection, penalization, backward stochastic differential
equation, Bellman-Isaacs equation, Perron's method.


\section{Introduction}

In this paper we study systems of variational inequalities with interconnected lower and upper obstacles. This type of inequalities arises as the Bellman-Isaacs equation in a multiple modes switching game between two players. Besides their classical fields of applications, multiple modes switching games are attracting a lot of interest in the management of power plants (see Bernhart (2011) (\cite{marieb}) and Perninge (2011) (\cite{perninge}), where they are successfully used to design optimal stopping and starting strategies for power flow control through activation of regulating bids on a regulated power market.

\ms\no The objective of this work is to establish existence and
uniqueness of a continuous viscosity solution of the following
system of variational inequalities with oblique reflection: \be
\label{mainsyst_varineq}  \left\{
\begin{array}{l}
\min\left\{ \big(v^{ij}- L^{ij}[\vec{v}]\big)(t,x),\,
\right.\max\left\{ \big(v^{ij} - U^{ij}[\vec{v}]\big)(t,x) ,\right.\\
\left.\left.\;\qquad \qq -\partial_tv^{i,j}(t,x)- {\cal
L}v^{ij}(t,x)-f^{ij}(t,x,(v^{kl}(t,x))_{(k,l)\in \Gamma^{1}
\times \Gamma^{2}}, \sigma^{\top}(t,x)D_x
v^{ij}(t,x))\right\}\right\} =0,\\
v^{ij}(T,x)=h^{ij}(x),
\end{array}\right. \ee
for every pair $(i,j)$ in the finite set of modes $\Gamma^{1}\times \Gamma^{2}$, where, for any $(t,x)\in \esp$,
$$
\dis{\mathcal{L}\varphi(t,x) := b(t,x) D_{x}\varphi(t,x) + \frac{1}{2}\textrm{Tr}[\sigma
\sigma^{\top}(t,x)D_{xx}^{2}\varphi(t,x)]},
$$
and to any solution $ \vec{v} =(v^{ij})_{\ij}$ we associate the
lower obstacle operator
$$
 \displaystyle{L^{ij}[\vec{v}] (t,x):= \max_{k \in\Gamma^{1},\,k\neq i} \big(v^{kj} -\underline{g}_{ik}\big)(t,x)},
 $$
and the upper obstacle operator
$$
\displaystyle{U^{ij}[\vec{v}](t,x) := \min_{l \in\Gamma^{2},\,l\neq j} \big(v^{il} +\bar{g}_{jl}\big)(t,x)},
$$
where, $\underline{g}_{ik}$ (resp. $\bar{g}_{jl}$) stands for the
switching cost incurred when the first (resp. second) player decides
to switch from mode $i$ to mode $k$ (resp. from mode $j$ to mode
$l$). Finally, the function $f^{ij}$ stands for the instantaneous payoff when the first
player is in mode $i$ and the second one in mode $j$.
\ms

\no The system (\ref{mainsyst_varineq}) and related switching games
have been studied by several authors. The most recent work
discussing this topic includes the papers by Hu and Tang (2008)
(\cite{Hutang-2}) and Tang and Hou (2007) (\cite{TangHou}) (see also
the references therein) which deal with switching games related to
(\ref{mainsyst_varineq}), when the switching costs do not depend on
the state variable. To the best of our knowledge, Ishii and Koike
(1991) (\cite{Ishiikoike91}) are the latest most general existence
and uniqueness results for the system (\ref{mainsyst_varineq}), for
state-dependent switching costs. They derive existence and
uniqueness of viscosity solutions of the elliptic version of
(\ref{mainsyst_varineq}) in a bounded domain of $\R^k$ whose
boundary is of class ${\cal C}^2$,  when the so-called {\it Fichera
functions} are strictly negative (see \cite{Ishiikoike91}, Proposition
4.3).\ms

\no The main result of the present paper, which is given in Theorem
\ref{mainthm2}, establishes existence and uniqueness of a continuous
viscosity solution of the system (\ref{mainsyst_varineq}), when the
state space  is the whole $\R^k$ and under rather weak assumptions
on the involved coefficients. Our approach is probabilistic and
makes use of penalization schemes that allow us to connect the
related penalized PDEs with systems of reflected backward SDEs with
unilateral interconnected obstacles which, for instance, have been
studied  in \cite{djehicheetal09}, \cite{Ham_elsasri09},
\cite{Hamzhang} or \cite{Hutang}. With the help of these sequences
of solutions of reflected BSDEs and their connection with PDEs, via
Feynman-Kac's formula, we are able to construct in Propositions
\ref{main-subsol} and \ref{main-supersol} both a viscosity subsolution
and a supersolution for the system (\ref{mainsyst_varineq}). Relying next both on
the comparison result established in Theorem \ref{comparison} and adapting the
Perron's method we construct a solution for
(\ref{mainsyst_varineq}) which is therefore unique. Finally, using again the uniqueness result, we identify the limit of the
penalized decreasing scheme as the solution of the system
(\ref{mainsyst_varineq}).\ms

\no We made this detour instead of trying to solve a related  system
of reflected BSDEs with interconnected bilateral obstacles, as one
would expect, simply because they satisfy neither the so-called
Mokobodski condition nor the condition of strict separation of the
two obstacles, which would guarantee existence and uniqueness of the
solutions of the system of BSDEs, since these obstacles  depend on
the solution. The structure of these bilateral obstacles suggests a
rather new type of conditions to guarantee existence and uniqueness
result for the related system of reflected BSDEs. This problem is
beyond the scope of the present paper and therefore left for future
research.

\ms\no Our plan for this paper is as follows. In Section 2, we
provide all the notations used in the paper, state the whole list of
required assumptions and define viscosity sub- and supersolutions along with
equivalent characterizations. In Section 3 we construct two
approximation schemes (an increasing and a decreasing on), consisting of
sequences of penalized reflected BSDEs associated with standard
switching problems. The counterpart of the decreasing 
scheme  (resp. the increasing one) in terms of PDEs stands 
for the penalized scheme of system (\ref{mainsyst_varineq}) (resp. system (\ref{approxdessous-1}) given below). 
Section 4 is devoted to the proof of a comparison result
related to the sub- and supersolutions of (\ref{mainsyst_varineq}).
In Section 5, the decreasing limit is identified as a viscosity
subsolution  of (\ref{mainsyst_varineq}) while a super-solution is
exhibited. Finally, we use Perron's method to construct a viscosity
solution of (\ref{mainsyst_varineq}) and,  thanks to the uniqueness
result, its connection with the limit of the decreasing scheme is
obtained. As a by product, we also obtain a similar characterization of the limit of the increasing scheme as
the unique solution in viscosity sense of the following system of variational inequalities: For every $\ij$, \be
\lb{approxdessous-1}\left\{
\begin{array}{l}
\max\left\{ \big(\underline v^{ij}-U^{ij}[\underline{\vec{v}}]\big)(t,x),\;\right.\min \left\{ \big(\underline
v^{ij}(t,x)-L^{ij}[\underline{\vec{v}}]\big)(t,x), \right. \\ \left. \left.\qquad \qq -\partial_t
\underline v^{ij}(t,x)- {\cal L}\underline
v^{ij}(t,x)-f^{ij}(t,x,(\underline v^{kl}(t,x))_{(k,l)\in
\Gamma^{1} \times \Gamma^{2}}, \sigma^{\top}(t,x)D_x
\underline v^{ij}(t,x))\right\}\right\} =0,\\
\underline v^{ij}(T,x)=h^{ij}(x).
\end{array}\right.\ee
We do not know whether the solutions of (\ref{mainsyst_varineq}) and (\ref{approxdessous-1}) coincide or not. We note that this is a very important 
issue since this will enable to characterize this common solution as the value of the zero-sum switching game. Since it is beyond the scope of the paper, this question 
is left for further research.
\section{Preliminaries and notation}
 Let $T$ (resp. $k$, $d$) be a fixed positive constant (resp.
two integers) and $\Gamma^{1}$ (resp. $\Gamma^{2} $) denote the set
of switching modes for player 1 (resp. 2). For later use, we shall
denote by $\Lambda$ the cardinal of the product set
$\Gamma^{1}\times \Gamma^{2} $ and for $(i,j)\in \gam$,
$\gami:=\Gamma^1-\{i\}$ and $\gamj:=\Gamma^2-\{j\}$. Next, for
$\vec{y}=(y^{kl})_{(k,l)\in \gam} \in \R^{\L}$, $\ij$, and
$\underline{y}\in \R$, we denote by $[\vec{y}^{-(ij)},\underline{y}]$
the matrix which is obtained from $\vec{y}$ by replacing the element
$y^{ij}$ with $\underline{y}$.

\ms\no Next, let us introduce the following functions.
 For any $\ij$,
$$
\begin{array}{l}
b:\, (t,x)\in [0,T]\times \R^k\mapsto b(t,x)\in \R^{k};\\
\sigma: \, (t,x)\in [0,T]\times \R^k\mapsto \sigma(t,x)\in \R^{k\times
d};\\ f^{ij}: \, (t,x,\vec{y},z)\in [0,T]\times
\R^{k+\Lambda+d}\mapsto
f^{ij}(t,x,\vec{y},z)\,\,\in \R\,;\\
\underline{g}_{ik}: \, (t,x)\in [0,T]\times\R^k\mapsto \underline{g}_{ik}(t,x)\in \R\,\;\quad (k \in (\Gamma^{1})^{-i});\\
\bar{g}_{jl}: \, (t,x)\in [0,T]\times\R^k \mapsto \bar{g}_{jl}(t,x)\in \R\,\;\quad (l \in (\Gamma^{2})^{-j});\\
h^{ij}: \, x\in \R^k\mapsto  h^{ij}(x)\in \R.
\end{array}
$$
\noindent A function $\Phi: \, (t,x)\in \esp\mapsto
\Phi(t,x)\in \R$ is called of {\it polynomial growth} if there exist
two non-negative real constants $C$ and $\g$ such that
$$ |\Phi(t,x)|\leq
C(1+|x|^\g),\quad (t,x)\in \esp.
$$
Hereafter, this class of functions is denoted by
$\Pi_{g}$. Let $\cC^{1,2}(\esp)$ (or simply
$\cC^{1,2}$) denote the set of real-valued functions defined on $\esp$,
which are once (resp. twice) differentiable w.r.t. $t$ (resp. $x$)
and with continuous derivatives.

\medskip\noindent In this paper, we  investigate
existence and uniqueness of viscosity solutions
$\vec{v}(t,x):=(v^{kl}(t,x))_{(k,l)\in \Gamma^{1} \times
\Gamma^{2}}$ of the following system of variational inequalities
with upper and lower interconnected obstacles: For any $\ij$,
\be \label{main-systvi}  \left\{
\begin{array}{l}
\min\left\{ \big(v^{ij}- L^{ij }[\vec{v}]\big)(t,x),\;
\right.\max\left\{ \big(v^{ij} - U^{ij}[\vec{v}]\big)(t,x),\right.\\
\left.\left.\;\qquad \qq -\partial_tv^{ij}(t,x)- {\cal
L}v^{ij}(t,x)-f^{ij}(t,x,(v^{kl}(t,x))_{(k,l)\in \Gamma^{1}
\times \Gamma^{2}}, \sigma^{\top}(t,x)D_x
v^{ij}(t,x))\right\}\right\} =0\\
v^{ij}(T,x)=h^{ij}(x)
\end{array}\right.
\ee
where, for any $(t,x)\in \esp$,
$$
\dis{\mathcal{L}\varphi(t,x) := b(t,x) D_{x}\varphi(t,x) + \frac{1}{2}\textrm{Tr}[\sigma
\sigma^{\top}(t,x)D_{xx}^{2}\varphi(t,x)] },
$$
and $\ij $,
$$
 L^{ij}[\vec{v}](t,x):=\max\limits_{k\in (\Gamma^{1})^{-i}}(v^{kj}(t,x)-\underline{g}_{ik}(t,x) ) \,\,
 \mbox{ and }\,\, U^{ij}[\vec{v}](t,x)= \min\limits_{l\in (\Gamma^{2})^{-j}}(v^{il}(t,x) +\bar{g}_{jl}(t,x)).
$$
The functions $f^{ij}$ stand for the instantaneous payoff when the first
player is in mode $i$ and the second one in mode $j$, and $
\underline{g}_{ik}$ (resp. $\bar{g}_{jl}$) stands for the
switching cost incurred when the first (resp. second) player decides
to switch from mode $i$ to mode $k$ (resp. from mode $j$ to mode
$l$).

\ms\no
The lower obstacle $L^{ij}[\vec{v}]$ and an upper obstacle $U^{ij}[\vec{v}]$ are called interconnected because each of them depends on the underlying solution $\vec{v}:=(v^{kl})_{(k,l)\in \Gamma^{1} \times\Gamma^{2}}$.

\medskip\noindent In a way, the system (\ref{main-systvi}) is the Bellman-Isaacs system of equations
associated with the zero-sum switching game with utility functions
$(f^{ij})_{\ij}$, terminal payoffs $(h^{ij})_{\ij}$ and switching
costs for the maximizer (resp. minimizer) given by $(\underline
g_{ij})_{\ij})$ (resp. $(\bar g_{ij})_{\ij})$). \ms

\medskip
The following assumptions are in force throughout the rest of the paper.
  \begin{enumerate}
\item[$\mathbf{(H0)}$] The functions $b$ and $\sigma$ associated with the second order operator $\mathcal{L}$
are jointly continuous in $(t,x)$, of linear growth in $(t,x)$ and
Lipschitz continuous w.r.t. $x$, meaning that there exists a
non-negative constant $C$ such that for any $(t,x,x')
 \in  [0,T] \times \mathbb{R}^{k+k}$ we have:
$$ |b(t,x) |+|\sigma(t,x)| \le C(1 +|x|) \;\,\,\textrm{and} \;\,\, |\sigma(t,x)- \sigma(t,x')|+ |b(t,x)- b(t,x')|
\le C|x -x'|.
$$
\item[$\mathbf{(H1)}$] Each function $f^{ij}$

(i) is continuous in $(t,x)$ uniformly w.r.t. the other variables
$(\vec{y},z)$ and for any $(t, x)$ and the mapping $(t,x)
\rightarrow f^{ij}(t,x,0,0)$ is of polynomial growth.

(ii) satisfies the standard hypothesis of Lipschitz continuity with
respect to the variables ($\vec{y}:=(y^{ij})_{(i,j)\in
\Gamma_1\times \Gamma_2},z$), i.e.
 $\forall \; (t,x) \in [0,T] \times \mathbb{R}^{k},\; \forall \; (\vec{y}_{1}, \vec{y}_{2}) \in \mathbb{R}^{\Lambda} \times  \mathbb{R}^{\Lambda}
, (z^{1}, z^{2}) \in \mathbb{R}^{d} \times \mathbb{R}^{d},$
$$ |f^{ij}(t,x, \vec{y}_{1},z_{1}) -f^{ij}(t,x, \vec{y}_{2},z_{2}) | \le C\left( |\vec{y}_{1}-\vec{y}_{2}| +|z_{1}-z_{2}|\right),$$
where, $|\vec{y}|$ stands for the standard Euclidean norm of
$\vec{y}$ in $\R^\L$.
\item[$\mathbf{(H2)}$] \underline{Monotonicity}: Let $\vec{y} = (y^{kl})_{(k,l) \in \Gamma^{1}\times \Gamma^{2}} $,
then for any $\ij$ and any $(k,l) \neq (i,j)$ the mapping $y^{k,l}
\rightarrow f^{ij}(s,\vec{y}, z)$ is non-decreasing.
\item[$\mathbf{(H3)}$] The functions $h^{ij}$, which are the terminal conditions in the system (\ref{main-systvi}), are continuous with respect to $x$, belong to class $ \Pi_{g}$ and satisfy
$$ \displaystyle{ \forall \; (i,j) \in \Gamma^{1} \times \Gamma^{2} \mbox{ and }x\in \R^k, \,\, \max_{k\in
(\Gamma^{1})^{-i}}\big( h^{kj}(x) -\underline{g}_{ik}(T,x) \big)
\le h^{ij}(x)  \le \min_{l\in (\Gamma^{2})^{-j}} \big(h^{il}(x)
+\bar{g}_{jl}(T,x)\big).} $$
\item[$\mathbf{(H4)}$] \underline{The no free loop property}: The switching costs $\underline{g}_{ik} $ and
$ \bar{g}_{jl} $ are non-negative, jointly continuous in $(t,x)$,
belong to $\Pi_g$ and satisfy the following condition:

For any loop in $\gam$, i.e., any sequence of pairs $(i_1,j_1),\ldots,(i_N,j_N)$
 of $\gam$ such that $(i_N,j_N)=(i_1,j_1)$, \mbox{card}$\{
(i_1,j_1),\ldots,(i_N,j_N)\}=N-1$ and $\forall \,\,q=1,\ldots,N-1$, either
$i_{q+1}=i_q$ or $j_{q+1}=j_q$, we have $\forall (t,x)\in \esp$,
\be\label{nonfreeloop3} \sum_{q=1,N-1}\varphi_{i_qi_{q+1}}(t,x)\neq
0, \ee where, $\forall \,\,\,q=1,\ldots,N-1,\,\,
\varphi_{i_qi_{q+1}}(t,x)=-\underline{g}_{i_qi_{q+1}}(t,x)\ind_{i_q\neq
i_{q+1}}+\bar{g}_{j_qi_{q+1}}(t,x)\ind_{j_q\neq j_{q+1}}$ \\(resp.
$
\varphi_{i_qi_{q+1}}(t,x)=\underline{g}_{i_q,i_{q+1}}(t,x)\ind_{i_q\neq
i_{q+1}}-\bar{g}_{j_q,i_{q+1}}(t,x)\ind_{j_q\neq j_{q+1}}).$

\end{enumerate}

This assumption implies in particular that \be \label{lp1}
\displaystyle{\forall \; (i_{1}, \ldots, i_{N})\in (\Gamma^1)^N
\;\textrm{such that } \;i_{N} =i_{1}\mbox{ and } \mbox{card}\{i_{1},
\ldots, i_{N}\}=N-1,\,\, \sum_{p=1}^{N-1}\underline{g}_{i_{k},
i_{k+1}}
> 0 } \ee and \be \label{lp2} \displaystyle{\forall \; (j_{1},
\ldots, j_{N})\in (\Gamma^2)^N \;\textrm{such that } j_{N} =j_{1}
\mbox{ and } \mbox{card}\{j_{1}, \ldots, j_{N}\}=N-1, \;
\sum_{p=1}^{N-1}\bar{g}_{j_{k}, j_{k+1}} > 0}.\ee By convention we
set $\bar{g}_{j,j} = \underline{g}_{i,i} =0$.

\ms\no Conditions (\ref{lp1})
and (\ref{lp2}) are classical in the literature of switching
problems and usually referred to as the {\it no free loop property}.

\ms\no Finally, let us mention that if we set
$$
\underline{g}_{ij}(t,x)=|i-j|\underline{g}(t,x)\;\;\,\mbox{and} \;\;\, \bar
g_{ij}(t,x)=|i-j|\bar{g}(t,x), \q (i,j)\in \gam, $$ where both
$\underline{g}$ and $\bar{g}$ are functions such that, for any
$\tx$, $\frac{\bar{g}(t,x)}{\underline{g}(t,x)}$ is not a rational
number, then assumption (\ref{nonfreeloop3}) holds.

\ms\no We now define the notions of viscosity super (or sub)-solution of the
system (\ref{main-systvi}). This is done in terms of the
notions of {\it subjet} and {\it superjet} which we recall here.

\begin{axiom}(Subjet and superjet) \\ \ms
$(i)$ For a lower semicontinuous (lsc) (resp. upper semicontinuous
(usc)) function $u:\, [0,T]\times \R^k\rightarrow \R,$ we denote by
$J^-u(t,x)$ (resp. $J^+u(t,x)$) the parabolic subjet (resp.
superjet) of $u$ at $(t,x)\in [0,T]\times \R^k$, as the set of
triples $(p,q,M)\in \R\times \R^k\times \mathbb{S}^k$ satisfying
\beaa &u(t',x')\geq \,(\mbox{resp.}
\leq)\,\,u(t,x)+p(t'-t)+q^\top(x'-x)+\\&\qq \qq\qq\qq \qq\qq
\qq\frac{1}{2}(x'-x)^\top M(x'-x)+o\big(|t'-t|+|x'-x|^2\big)\eeaa
where $\mathbb{S}^k$ is the set of symmetric real matrices of
dimension $k$. \bs

\no$(ii)$ For a lsc (rep. usc) function $u:[0,T]\times \R^k\rightarrow \R$, we denote by $\bar J^-u(t,x)$ (resp. $\bar J^+u(t,x)$)
the parabolic limiting subjet (resp. superjet) of $u$ at $(t,x)\in
[0,T]\times \R^k$, as the set of triples $(p,q,M)\in \R\times
\R^k\times \mathbb{S}^k$ such that:$$ \ba{l} (p,q,M)=\lim_n
(p_n,q_n,M_n), \,\, (t,x)=\lim_n (t_n,x_n) \mbox{ with
}(p_n,q_n,M_n)\in J^{-}u(t_n,x_n)\\ \,\,(\mbox{resp. }J^+u(t_n,x_n))
\mbox{ and }u(t,x)=\lim_{n} u(t_n,x_n).  \ea$$
\end{axiom}

\medskip\noindent Finally, given a locally bounded
$\R$-valued deterministic function $u$, we denote by $u_{*}$ (resp.
$u^{*}$) its
 lower (resp. upper) semicontinuous envelope defined by:
\begin{equation}\label{eq:definitions} \forall \; (t,x)\in \esp,  \quad \dis{u_{*}(t,x) =
\liminf_{(t^{'}, x^{'}) \to (t,x);\,\, t'<T} u(t', x') \;\,\,\textrm{and} \;\,\,
u^{*}(t,x) = \limsup_{(t', x') \to (t,x);\,\, t'<T} u(t', x')}.
\end{equation}

\ms\no We now give the definition of a viscosity solution for the system (\ref{main-systvi}).

\begin{axiom} Viscosity solution to (\ref{main-systvi})\\
(i) A function $\vec{v} =(v^{kl}(t,x))_{(k,l)\in \gam}:[0,T]\times
\R^{k}\mapsto \R^{\Lambda}$ such that for any $(i, j)\in \Gamma^{1}
\times \Gamma^{2}$, $v^{ij}$ is lsc (resp. usc), is called a
viscosity supersolution (resp. a viscosity subsolution) to
(\ref{main-systvi}) if for any $(i,j)\in \Gamma^{1} \times
\Gamma^{2}$, for any $(t,x)\in [0,T)\times \R^{k}$ and any
$(p,q,M)\in \bar J^{-}v^{ij}(t,x)$ (resp. $\bar J^{+}v^{ij}(t,x)$)
we have: \be\label{defsurso}\left\{\begin{array}{l}
\min\left\{v^{ij}(t,x)-L^{ij}[\vec{v}](t,x),\right.\\
\left.\quad \max\left\{-p-b(t,x).q-\frac{1}{2}
Tr[(\sigma\sigma^\top)(t,x)M]-f^{ij}(t,x,\vec{v}(t,x),\sigma^\top(t,x)q)
; v^{ij}(t,x)-U^{ij}[\vec{v}](t,x)\right\} \right\}\geq 0,
\\v_i(T,x)\geq \, \, h^{ij}(x),\end{array}\right.\ee \Big(\mbox{resp.} \be
\label{defsousso}\left\{\begin{array}{l}
\min\left\{v^{ij}(t,x)-L^{ij}[\vec{v}](t,x),\right.\\
\left.\quad \max\left \{-p-b(t,x).q-\frac{1}{2}
Tr[(\sigma\sigma^\top)(t,x)M]-f^{ij}(t,x,\vec{v}(t,x),\sigma^\top(t,x)q)
; v^{ij}(t,x)-U^{ij}[\vec{v}](t,x)\right\}\right\}\leq 0,
\\ v^{ij}(T,x)\leq \; h^{ij}(x)\Big).
\end{array}\right.\ee
 (ii) A locally bounded function $\vec{v} = (v^{kl})_{(k,l) \in \Gamma^{1}
 \times \Gamma^{2}}:[0,T]\times \R^k\rightarrow \R^{\Lambda}$ is
called a viscosity solution of (\ref{main-systvi}) if the associated
lower (resp. upper) semicontinuous envelope
${(v^{kl}_{*})}_{(k,l)\in \gam}$ (resp. $({v^{kl}}^{*})_{(k,l)\in
\gam}$) defined in (\ref{eq:definitions}) is a viscosity
supersolution (resp. subsolution) of (\ref{main-systvi}).

\noindent If, in addition, for any $(k,l)\in \gam$, $ v^{kl}_{*}
=v^{kl *}$, then $(v^{kl})_{(k,l)\in \gam}$ is a continuous
viscosity solution of (\ref{main-systvi}).
\end{axiom}
\begin{remark}\lb{allegement} Under Assumptions (H0)-(H4), the above definition of a viscosity
solution for (\ref{main-systvi}) can be relaxed replacing the limiting subjet $\bar{J}^{-}(v_{*}^{ij})(t,x) $ of the supersolution
 $v_{*}$ (resp. the limiting superjet $\bar{J}^{+}(v^{ij, *})(t,x) $ of $v^{*}$) by the subjet $J^{-}(v_{*}^{ij})(t,x)$ (resp. by the superjet
 $J^{+}(v^{ij, *})(t,x)$). This results from the continuity of the functions $b$, $\sigma$, $(f^{ij},
h^{ij},\underline g_{ij}, \bar g_{ij})_{\ij}$ and the monotonicity
property (H2) of $(f^{ij})_{\ij}$. 
\end{remark}


\section{Systems of reflected BSDEs and approximation schemes of the solutions }

Let $(\Omega, {\cal F}, \P)$ be a fixed probability space on which
is defined a standard $d$-dimensional Brownian motion
$B=(B_t)_{0\leq t\leq T}$ whose natural filtration is
$(\cF_t^0:=\sigma \{B_s, s\leq t\})_{0\leq t\leq T}$. Let $
\bF=(\cF_t)_{0\leq t\leq T}$ be the completed filtration of
$(\cF_t^0)_{0\leq t\leq T}$ with the $\mathbb{P}$-null sets of
${\cal F}$, hence $(\cF_t)_{0\leq t\leq T}$ satisfies the usual
conditions, i.e., it is right continuous and complete.
Furthermore, let
\begin{itemize}
\item ${\cal P}$ be the $\sigma$-algebra on $[0,T]\times \Omega$ of
$\bF$-progressively measurable sets;

\item ${\cal H}^{2,\ell}$ ($\ell \geq 1$) be the set of $\cal
P$-measurable and $\R^\ell$-valued processes $w=(w_t)_{t\leq T}$
such that $\E[\int_0^T|w_s|^2ds]<\infty$;

\item ${\cal S}^{2,\ell}$ ($\ell \geq 1$) be the subset of ${\cal
H}^{2,\ell}$ of continuous processes such that $\E[\sup_{t\leq
T}|{w}_t|^2]<\infty$. Finally let $\mathcal{A}^{+,2}$ be the subset
of ${\cal S}^{2,1}$ of non-decreasing processes $K=(K_t)_{t\leq T}$
such that $K_0=0$ and $\E[K_T^2]<\infty$.
\end{itemize}
\medskip

\noindent Next, for $n,m\geq 0$, let $(Y^{ij,n,m},Z^{ij,n,m})_{(i,j)\in \gam}$
be the solution of the following system of BSDEs.
 \be
\label{doublypenalizedscheme} \left\{
\begin{array}{l}
(Y^{ij,n,m},Z^{ij,n,m})\in {\cal
S}^{2,1}\times {\cal H}^{2,d};\\
 dY_{s}^{ij,n,m} \; = \; -f^{ij,n,m}(s,X^{t,x}_s, (Y_{s}^{kl,n,m})_{(k,l)\in \gam},
 Z_{s}^{ij,n,m}) ds + Z_{s}^{ij,n,m }dB_{s},\,\,\,  s\leq T,\\
Y_{T}^{ij,nm} = h^{ij}(X^{t,x}_{T}),
\end{array} \right.
\ee where, $$\ba{l} f^{ij,n,m}(s, X^{t,x}_s, (y^{ij})_{(i,j)\in
\gam},z):=f^{ij}(s, X^{t,x}_s, (y^{kl})_{(k,l)\in
\gam},z)+\\\qq\qq n(y^{ij}-\max_{k\in \gami
}\{y^{kj}-\underline{g}_{ik}(s,X^{t,x}_s)\})^{\mathbf{-}}-m(y^{ij}-\min_{l\in
\gamj }\{y^{il}+\bar{g}_{jl}(s,X^{t,x}_s)\})^{\mathbf{+}}. \ea$$
Note that by Assumption (H1) and the standard result on
multi-dimensional BSDEs by Pardoux and Peng (1990) (\cite{Pardouxpeng}), the
solution $(Y^{ij,n,m},Z^{ij,n,m})_{(i,j)\in \gam}$ exists and is
unique. Since the generators and the terminal values depend on
$(t,x)$, the processes $Y^{ij,n,m}$ and $Z^{ij,n,m}$ also depend on $(t,x)$
but, to avoid overload notation, we do not mention this dependence in the sequel.
Furthermore,  the following monotonicity
 properties holds for the double sequence $(Y^{ij,n,m})_{n,m}$.
\begin{propo} \label{comparaison} For any $\ij$ and $n,m\geq 0$ we have
\be\lb{ing}
\P-a.s.,\,\,\quad Y^{ij,n,m}\leq Y^{ij,n+1,m}\,\,\, \mbox{ and }\,\,\, Y^{ij,n,m}\geq
Y^{ij,n,m+1},\qquad (i,j)\in \gam.\ee
Moreover, for any $(i,j)\in \gam$ and $n,m\geq 0$, there exists a deterministic
continuous function $v^{ij,n,m}$ in $\Pi_g$ such that, for any $t\leq T$,
\be \lb{ynm}Y^{ij,n,m}_s=v^{ij,n,m}(s,X^{t,x}_s),\quad s\in [t,T].
\ee
Finally, for any $(i,j)\in \gam$ and $n,m\geq 0$,
\be \lb{ineqvs}v^{ij,n,m}(t,x)\le
v^{ij,n+1,m}(t,x) \,\,\, \mbox{ and }\,\,\, v^{ij,n,m}(t,x)\ge v^{ij,n,m+1}(t,x),\qquad (t,x)\in [0,T]\times \R^k.\ee
\end{propo}
\pr.
The second claim is just the representation of solutions of standard
BSDEs by deterministic functions in the Markovian framework (see e.g. El Karoui {\it et al.} (1997) (\cite{elkarouiquenezpeng}) for more details). As
for the first one, it is based on the result by Hu and Peng (2006)
(\cite{Hupeng}) related to the comparison of solutions of
multi-dimensional BSDEs (we recall in  Appendix (A1)).
To prove this, it is enough to show that for any $t$,
$({y}_{ij})_{(i,j)\in \gam}$, $({\bar y}_{ij})_{(i,j)\in \gam}$ $\in
\R^{\Lambda}$ and $(z_{ij})_{(i,j)\in \gam}$, $({\bar
z}_{ij})_{(i,j)\in \gam}$ $\in (\R^d)^\Lambda$ we have:
$$\ba{l}
-4\sum_{(i,j)\in \gam}y_{ij}^-\left (f^{ij,n,m}(s, X^{t,x}_s,
(y_{kl}^+ +{\bar y}_{kl})_{(k,l)\in \gam},z_{ij})-f^{ij,n+1,m}(s,
X^{t,x}_s, ({\bar y}_{kl})_{(k,l)\in \gam},\bar z_{ij})\right)\\
\qq\qq\qq \leq 2\sum_{(i,j)\in \gam}1_{[y_{ij}<0]}\|z_{ij}-\bar
z_{ij}\|^2+C\sum_{(i,j)\in \gam}(y_{ij}^-)^2, \ea
$$
where, $C$ is a constant, $y_{ij}^-=\max(-y_{ij},0)$ and
$y_{ij}^+=\max(y_{ij},0)$.  This
inequality follows easily from the fact that, for any $(i,j)\in \gam$,

$(i)$ $f^{ij,n,m}(s, X^{t,x}_s, ({y}_{kl})_{(k,l)\in
\gam},z_{ij})\leq f^{ij,n+1,m}(s, X^{t,x}_s, ({\bar
y}_{kl})_{(k,l)\in \gam},z_{ij}),$

$(ii)$ For any $(u_{kl})_{(k,l)\in \gam}\in \R^\Lambda$ such that
$u_{ij}=0$ and $u_{kl}\geq 0$, for $(k,l)\neq (i,j)$,
$$f^{ij,n,m}(s, X^{t,x}_s, ({y}_{kl}+u_{kl})_{(i,j)\in
\gam},z_{ij})\geq f^{ij,n,m}(s, X^{t,x}_s, ({ y}_{kl})_{(k,l)\in
\gam},z_{ij}).$$

$(iii)$ $f_{ij}$ depends only on $z_{ij}$ and not on the other
components $z_{kl}$, $(k,l)\neq (i,j)$.

\noindent Consequently, we have
$$
\P-a.s., \qquad Y^{ij,n,m}\leq Y^{ij,n+1,m},\qquad (i,j)\in \gam.
$$ In the same way one can show that $
Y^{ij, n,m+1}\leq Y^{ij,n,m}$.  Finally, the inequalities of (\ref{ineqvs})
are obtained by taking $s=t$ in (\ref{ing}) in view of the
representation (\ref{ynm}) of $Y^{ij,n,m}$ by $v^{ij, n,m}$ and $\x$.\qed

\ms We now suggest two approximation schemes obtained from the sequence $(Y^{ij,n,m},\,\,\ij)_{n,m}$ of the solution of the system (\ref{doublypenalizedscheme}). The first scheme is a sequence of decreasing reflected BSDEs with interconnected lower obstacles and the second one is an increasing sequence of reflected BSDEs with interconnected upper obstacles.

\ms Let us first introduce the decreasing approximation scheme by considering the following
system of reflected BSDEs with interconnected obstacles: $\forall
(i,j)\in \gam$, \be\label{multidimsystem}\left\{\begin{array}{l}
\bar Y^{ij,m}\in \cS^{2,1}, \,\,\bar Z^{ij,m} \in \cH^{2,d} \mbox{
and
}\bar K^{ij,m } \in \cA^{2,+}\,\,;\\
\bar Y^{ij,m}_s=h^{ij}(X^{t,x}_T)+\int_{s}^{T}\bar
{f}^{ij,m}(r,X^{t,x}_r,
(\bar Y^{kl,m}_r)_{(k,l)\in \gam},\bar Z^{ij,m}_r)dr+\bar K^{ij,m}_T-\bar K^{ij,m}_s-\int_s^T\bar Z^{ij,m}_rdB_r,\\
\displaystyle{\bar Y^{ij,m}_{s}\geq \max_{k\in
(\Gamma^{1})^{-i}}\{\bar
Y^{kj,m}_s-\underline{g}_{ik}(s,X^{t,x}_s)\}},
\,\,\qquad s\leq T,\\
\int_0^T(\bar Y^{ij,m}_s-\max_{k\in (\Gamma^{1})^{-i}}\{\bar
Y^{kj,m}_s-\underline{g}_{ik}(s,X^{t,x}_s)\})d\bar
K^{ij,m}_s=0,\end{array}\right.\ee where, for any $(i,j) \in \gam$,
$m\geq 0$ and $(s, \vec{y}, z^{ij})$,
$$\bar{f}^{ij,m}(s, X^{t,x}_s,\vec{y},z^{ij}):= f^{ij}(s,X^{t,x}_s,(y^{kl})_{(k,l)\in \gam},z^{ij})
-m\big(y^{ij} -\min_{l\in \gamj} (y^{il} + \bar{g}_{jl}
(s,X^{t,x}_s))\big)^{+}.
$$
Thanks to Assumptions (H1)-(H3) and (\ref{lp1}), by  Theorem 3.5 in
Hamad\`ene and Zhang (2010) (\cite{Hamzhang}), the solution of
(\ref{multidimsystem}) exists and is unique. In fact, this holds again under
 weaker assumptions (see Hamad\`ene and Morlais (2011)(\cite{Hammorlais11}, 
Theorem 1). For sake of completeness, a statement of this recent result is given in Appendix (A2).  Moreover, we have the following properties.

\medskip

\begin{propo}\label{firstresult} For any $(i,j)\in \gam$ and $m\geq 0$, we have

\noindent (i)  \be \label{conv1} \E[\sup_{t\le s\leq
T}|Y^{ij,n,m}_s-\bar Y^{ij,m}_s|^2]\rw 0, \,\,\mbox{ as } \,\,{n\rw
\infty},\ee (ii)
$$
\P-a.s., \qquad \bar Y^{ij,m}\geq \bar Y^{ij,m+1}.
$$
(iii) There exists
a unique $\Lambda$-uplet of deterministic continuous functions
$(\bar v^{kl,m})_{(k,l)\in \gam}$ in $\Pi_g$ such that, for every $t\le T$,
\be\lb{exp} \bar
Y^{ij,m}_s=\bar v^{ij,m}(s,X^{t,x}_s),\qquad  s\in [t,T].\ee
Moreover, for any
$(i,j)\in \gam$ and $(t,x)\in [0,T]\times \R^k$, $\bar v^{ij,m}(t,x)\ge \bar v^{ij,m+1}(t,x)$.

\medskip\noindent
Finally,
$(\bar v^{ij,m})_{(i,j)\in \gam}$ is the unique viscosity solution in the class
$\Pi_g$ of the following system of
variational inequalities with inter-connected obstacles.
$\forall\,\, (i,j)\in \gam$,
\begin{equation}\label{newauxiliary-systvi}\left\{\ba{l}
\min\{\bar v^{ij,m} (t,x)- \max\limits_{k\in
(\Gamma^{1})^{-i}}(\bar
v^{kj,m}(t,x)-\underline{g}_{ik}(t,x));\qq\qq\qq
\\\qq-\partial_t
\bar v^{ij,m}(t,x)-b(t,x)D_x\bar
v^{ij,m}(t,x)-\frac{1}{2}\textrm{Tr}(\sigma \sigma^\top
(t,x)D^2_{xx} \bar v^{ij,m}(t,x))\\\qq\qq\qq\qq
-\bar{f}^{ij,m}(t,x,(\bar v^{kl,m}(t,x))_{(k,l)\in
\gam},\sigma^\top(t,x)D_x\bar v^{ij,m}(t,x))\} = 0,\\
\bar v^{ij,m}(T,x)=h^{ij}(x) \ea\right.
\end{equation} where,
$$
\bar{f}^{ij,m}(s, x, (y^{kl})_{(k,l)\in \gam},z) = f^{ij}
(s,x,(y^{kl})_{(k,l)\in \gam},z) -m\big(y^{ij} -\min_{l\in \gamj}
(y^{il} + \bar{g}_{jl}(s,x) )\big)^{+}.
$$
\end{propo}
$\pr.$ Let us prove (i). For this, it is enough to consider the case $m=0$,
and we will do so, since for any $\ij$, the function
$$(s,x,(y^{kl})_{(k,l)\in \gam}) \longrightarrow -m\big(y^{ij} -\min_{l\in
\gamj} (y^{il} + \bar{g}_{jl}(s,x)))^+$$ has the same properties
as $f^{ij}$ displayed in (H1)-(H2).

\medskip
\noindent To begin with, let us show that for any $i,j$ and $n\geq 0$,
\be\lb{ineqxx}
\P-a.s., \qquad Y^{ij,n,0}\le  \bar Y^{ij,0}.
\ee

\noindent  First andt w.l.o.g. we may assume that $f^{ij}$ is
non-decreasing w.r.t. $(y^{kl})_{(k,l)\in \gam}$, since thanks to Assumption (H2), it is enough to multiply the solutions by
$e^{\varpi t}$, where  $\varpi$ is appropriately chosen in order to fall
in this latter case, since $f^{ij}$ is Lipschitz in $y^{ij}$. Now,
for fixed $n$, let us define recursively the sequence $(\tilde
Y^{k,ij,n})_{k\ge 0}$ as follows:

\medskip\noindent For $k=0$ and any $\ij$, we set $
\tilde Y^{0,ij,n} :=\bar Y^{ij,0}$ and, for any $k\geq 1$, let us define $(\tilde Y^{k,ij,n},\tilde
Z^{k,ij,n})\in {\cal S}^{2,1}\times {\cal H}^{2,d}$ as the solution
of the following system of BSDEs: $\forall \ij$, \be
\label{doublypenalizedscheme3} \left\{
\begin{array}{l}
- d\tilde Y_{s}^{k,ij,n} \; = \; f^{ij}(s,X^{t,x}_s, (\tilde
Y_{s}^{k-1, pq,n})_{(p,q)\in \gam},
 \tilde  Z_{s}^{k,ij,n}) ds +\\\qq\qq\qq\qq n\{\tilde Y_{s}^{k,ij,n}-\max_{l\in \gami}(
\tilde Y_{s}^{k-1,lj,n}-\underline{g}_{ik}(s,X^{t,x}_s))\}^-ds
 - \tilde Z_{s}^{k,ij,n }dB_{s},\\
\tilde Y_{T}^{k,ij,n} = h^{ij}(X^{t,x}_{T})
\end{array} \right.
\ee
The solution of (\ref{doublypenalizedscheme3}) exists since it is a multi-dimensional 
standard BSDE with a Lipschitz coefficient, noting that $(\tilde Y_{s}^{k-1,
pq,n})_{(p,q)\in \gam}$ is already given. Since, $n$ is fixed and the coefficient
$$ \varphi^{ij,n}(s,\omega,(y^{kl})_{(k,l)\in \gam},z^{ij}):=f^{ij}(s,X^{t,x}_s,
(y^{kl})_{(k,l)\in \gam},z^{ij})+n\{y^{ij}-\max_{l\in \gami}(
y^{lj}-\underline{g}_{il}(s,X^{t,x}_s)\}^-
$$
is Lipschitz w.r.t. $((y^{kl})_{(k,l)\in \gam},z^{ij})$, the sequence $(\tilde Y^{k,ij,n})_k$
converges in ${\cal S}^{2,1}$ to $Y^{ij,n,0}$, as $k\to\infty$, for any $i,j$ and $n$.

\noindent Using an induction argument w.r.t. $k$, we prove that, for any $i,j$ and $n$,
$$
\P-a.s., \qquad\tilde Y^{k,ij,n}\leq \bar Y^{ij,0},\qquad k\ge 0.$$
Indeed,  for $k=0$ the
inequalities hold true through the definition of $\tilde
Y^{0,ij,n}$. Assume now that these inequalities are valid for some
$k-1$, i.e., for any $i,j$ and $n$,
\be \label{rechyp}
\P-a.s., \qquad \tilde Y^{k-1,ij,n}\leq \bar Y^{ij,0}.
\ee
Thus, taking into account that
$\bar Y^{ij,0}$ satisfies (\ref{multidimsystem}) with $m=0$ and the
fact that $f^{ij}$ is non-decreasing w.r.t. $(y^{kl})_{(k,l)\in
\gam}$ then for any $i,j,n$, it holds
$$
\ba{l}
f^{ij}(s,X^{t,x}_s, (\tilde Y_{s}^{k-1,pq,n})_{(p,q)\in
\gam},z^{ij})+n\{\bar Y^{ij,0}-\max_{l\in \gami}( \tilde
Y_{s}^{k-1,lj,n}-\underline{g}_{il}(s,X^{t,x}_s))\}^-\\\qq\qq \le
f^{ij}(s,X^{t,x}_s, (\bar Y_s^{pq,0})_{(p,q) \in \gam},z^{ij}). \ea
$$
Using the standard comparison result of solutions of one dimensional
BSDEs we obtain that, for any $n,i,j$, $\tilde Y^{k,ij,n}\leq \bar
Y^{ij,0}.$ a.s. Thus the property (\ref{rechyp}) is valid for any $k$. Taking the limit as $k$ tends to $\infty$, we obtain (\ref{ineqxx}). \ms

\ms\no We can now apply Peng's monotonic limit
 (see Peng (1999) (\cite{Peng99})) to the increasing sequence $(Y^{ij,n,0})_{n\ge 0}$. This yields
 the existence of processes $\hat Y^{ij}$, $\hat Z^{ij}$ and
$\hat K^{ij}$, $(i,j)\in \gam$, such that:

\noindent (a) $\hat Y^{ij}$ is RCLL and uniformly $\P$-square
integrable. Moreover, for any stopping time $\tau$, $\lim_{n\rw
\infty}\nearrow Y_\tau^{ij,n,0}=\hat Y^{ij}_\tau$.

\noindent (b) $\hat K^{ij}$ is RCLL non-decreasing,  $\hat
K^{ij}_0=0$ and for any stopping time $\tau$, $\lim_{n\rw
\infty}K_\tau^{ij,n,0}=\hat K^{ij}_\tau $, $\P-a.s.$

\noindent (c) $\hat Z^{ij}\in {\cal H}^{2,d}$ and for any $p
\in [1, 2)$,
$$
\lim_{n\rw \infty}\E[\int_0^T|Z_s^{ij,n,0}-\hat
Z^{ij}_s|^pds]=0.$$

\noindent (d) For any $(i,j)\in \gam$, the triple $(\hat Y^{ij}, \hat Z^{ij},\hat K^{ij})$ satisfies: $\forall
s\leq T$, \be\lb{eqint1}\left\{\ba{l} \hat Y_s^{ij}
=h^{ij}(X^{t,x}_T)+\int_{s}^{T} {f}^{ij}(r,X^{t,x}_r, (\hat
Y^{kl}_r)_{(k,l)\in \gam},\hat Z^{ij}_r)dr+\hat
K^{ij}_T-\hat K^{ij}_s- \int_s^T\hat
Z^{ij}_rdB_r\\ \hat Y^{ij}_{s}\geq \max_{k\in (\Gamma^{1})^{-i}}\{\hat
Y^{kj}_s-\underline{g}_{ik}(s,X^{t,x}_s)\} .\ea\right.\ee
The remaining of the proof is the same as in
Hamad\`ene and Zhang (2010) ((\cite{Hamzhang}), Theorem 3.2) and it mainly consists in proving both the continuity of $\hat Y^{ij}$ and
$\hat K^{ij}$ and the Skorohod condition 
(these properties are deduced using the no-free loop property (H4)). Finally we obtain that the triple $(\hat Y^{ij}, \hat Z^{ij},\hat K^{ij})$ satisfies
(\ref{multidimsystem}) with $m=0$. Hence, by uniqueness of the solutions of (\ref{multidimsystem}), $\hat Y^{ij}=\bar Y^{ij,0}$ a.s, which completes the proof of (i).

\medskip\noindent (ii) is an immediate consequence of (i) and Proposition
\ref{comparaison}. \ms

\medskip\noindent We now establish (iii). The existence of the deterministic continuous
functions $(\bar v^{ij,m})_{\ij}$ such that for any $s\in
[t,T],\,\,\, Y^{ij,m}_{s} = \bar v^{ij,m}(s, X_{s})$, $(i,j)\in
\gam$,  $m\ge 0$, and which satisfy the system
(\ref{newauxiliary-systvi}), is obtained in \cite{Hammorlais11} (see
Appendix A2, Theorem \ref{a2}). Finally as $Y^{ij,m}\geq
Y^{ij,m+1}$ a.s., we deduce that $\bar v^{ij,m}\geq \bar v^{ij,m+1}$,
taking into account (ii), which completes the proof. \qed \bs


We now introduce the increasing approximation scheme by considering the following system
of reflected BSDEs with interconnected obstacles: for any $(i,j)\in
\gam$, \be\label{multidimsystem2}\left\{\begin{array}{l} \underline
Y^{ij,n}\in \cS^2, \,\,\underline Z^{ij,n} \in \cH^{2,d} \mbox{ and
}\underline K^{ij,n } \in \cA^{2,+}\,\,;\\
\underline Y^{ij,n}_s=h^{ij}(X^{t,x}_T)+\int_{s}^{T}\underline
{f}^{ij,n}(r,X^{t,x}_r, (\underline Y^{kl,n}_r)_{(k,l)\in
\gam},\underline Z^{ij,n}_r)dr-(\underline K^{ij,n}_T-
\underline K^{ij,n}_s)-\int_s^T\underline Z^{ij,n}_rdB_r,\\
\displaystyle{\underline Y^{ij,n}_{s}\leq \min_{l\in
(\Gamma^{2})^{-j}}\{\underline
Y^{il,n}_s+\bar{g}_{jl}(s,X^{t,x}_s)\}},
\quad s\leq T, \\
\displaystyle{\int_0^T(\underline Y^{ij,n}_s-\min_{l\in
(\Gamma^{2})^{-j}}\{\underline
Y^{il,n}_s+\bar{g}_{jl}(s,X^{t,x}_s)\})d\underline
 K^{ij,n}_s=0,}\end{array}\right.
 \ee
where, for any $(i,j) \in \gam$, $n\geq 0$ and $(s, \vec{y}, z^{ij})$,
$$
\underline{f}^{ij,n}(s, X^{t,x}_s,\vec{y},z^{ij}) = f^{ij}(s,X^{t,x}_s,(y^{kl})_{(k,l)\in
\gam},z^{ij}) +n\big(y^{ij} -\max_{k\in \gami} \{y^{kj} -
\underline {g}_{ik}(s,X^{t,x}_s)\} \big)^{-}.
$$
Thanks to assumptions (H1)-(H3) and (H4)-(\ref{lp2}), by Theorems 3.2 and 3.5 in Hamad\`ene and Zhang (2010)  (\cite{Hamzhang})  the solution of
(\ref{multidimsystem}) exists and is unique.

\ms \noindent
Below is the analogous of Proposition \ref{firstresult}. We do not
give its proof since it is deduced from this latter proposition in
considering the equation satisfied by $(-\underline
Y^{ij,n},-\underline Z^{ij,n},\underline K^{ij,n})$.
\begin{propo}\label{secondresult} For any $(i,j)\in \gam$ we have

\noindent $(i)$ \be \label{conv2} \E[\sup_{t\le s\leq
T}|Y^{ij,n,m}_s-\underline{Y}^{ij,n}_s|^2]\rw 0 \quad \mbox{ as }\quad m\to
\infty.
\ee
$(ii)$ For any $n\geq 0$,
$$
\P-a.s., \qquad \underline{Y}^{ij,n}\leq \underline
{Y}^{ij,n+1}.
$$
$(iii)$ For every $n\ge 0$, there exists a unique
$\Lambda$-uplet of deterministic continuous functions $(\underline
v^{kl,n})_{(k,l)\in \gam}$ in $\Pi_g$ such that, for every $t\leq T$,
\be\lb{exp2} \underline
{Y}^{ij,n}_s=\underline v^{ij,n}(s,X^{t,x}_s),\qquad t\le s\le T.\ee
Moreover, for any
$n\geq 0$ and $(i,j)\in \gam$,
$$
\underline{v}^{ij,n}(t,x)\le
\underline{v}^{ij,n+1}(t,x),\qquad (t,x)\in[0,T]\times \R^k.
$$
Finally, $(\underline v^{ij,n})_{(i,j)\in \gam}$ is the unique
viscosity solution in the class $\Pi_g$ of the following system of
variational inequalities with inter-connected obstacles: $\forall
\,(i,j)\in \gam$,
\begin{equation}\label{newauxiliary-systvi-bis}\left\{\ba{l}
\max\{\underline v^{ij,n} (t,x)- \min\limits_{l\in
(\Gamma^{2})^{-j}}(\underline v^{il,n}(t,x)+\bar
{g}_{ik}(t,x));\qq\qq\qq
\\\qq-\partial_t
\underline v^{ij,n}(t,x)-b(t,x)D_x\underline
v^{ij,n}(t,x)-\frac{1}{2}\textrm{Tr}(\sigma \sigma^\top
(t,x)D^2_{xx} \underline v^{ij,n}(t,x))\\\qq\qq\qq\qq
-\underline{f}^{\,\,ij,n}(t,x,(\underline v^{kl,n}(t,x))_{(k,l)\in
\gam},\sigma^\top(t,x)D_x\underline v^{ij,n}(t,x))\} = 0,\\
\underline v^{ij,n}(T,x)=h^{ij}(x) \ea \right.
\end{equation}
where,
$$
\underline f^{\,\,\,ij,n}(s, x, (y^{kl})_{(k,l)\in \gam},z) = f^{ij}(s,x,(y^{kl})_{(k,l)\in \gam},z)
+n\{y^{ij} -\max_{k\in \gami} (y^{kj} - \underline{g}_{ik}(s,x)
)\}^{-}.
$$
\end{propo}

\medskip
For $\tx \in \esp$ and $\ij$, let us define
$$
\bar v^{ij}(t,x):=\lim_{m\rw \infty}\bar v^{ij,m}(t,x) \,\,\, \mbox{ and
}\,\,\, \underline v^{ij}(t,x):=\lim_{n\rw \infty}\underline v^{ij,n}(t,x).
$$
Then, as a by-product of Propositions \ref{firstresult} and
\ref{secondresult} we have

\begin{corollary}\lb{cndpol} For any $\ij$, the function $\bar v^{ij}$ (resp. $\underline v^{ij}$) is usc (resp. lsc). Moreover,  $\bar v^{ij}$ and $\underline v^{ij}$ belong to $\Pi_g$ i.e., there exist deterministic constants $C\geq 0$ and $\g>0$ such
that, for any $(t,x)\in \esp$,
$$
|\bar v^{ij}(t,x)|+|\underline
v^{ij}(t,x)|\leq C(1+|x|^\g),  \qquad(t,x)\in \esp.
$$

\ms \noindent Finally, for any $(t,x)\in \esp$,$$\underline v^{ij}(t,x)\leq \bar
v^{ij}(t,x).$$
\end{corollary}
$\pr$. For any $\ij$, $\bar v^{ij}$ (resp. $\underline{v}^{ij}$)
is obtained as a decreasing (resp. increasing) limit of continuous
functions. Therefore, it is usc (resp. lsc). Next, for any $(i,j)$
and $n,m$,
$$
v^{ij,n,m}(t,x)\le v^{ij,n,0}(t,x),\quad (t,x)\in [0,T]\times \R^k,
$$
as the sequence
$(v^{ij,n,m})_{m\ge 0}$ is decreasing. Thus, taking the limit as
$m\rw \infty$ we obtain
$$
\underline{v}^{ij,n}\le v^{ij,n,0}.
$$
Now, using  (\ref{ynm}) and (\ref{conv1}) it follows that, for any $t\le T$ and
$s\in [t,T]$, $Y^{ij,n,0}_s=v^{ij,n,0}(s,\x_s)$ and the processes
$Y^{ij,n,0}$ converges in $\cS^2$, as $n\rw \infty$, to $\bar
Y^{ij,0}$ solution of (\ref{multidimsystem}) with $m=0$. Furthermore,
by (\ref{exp}), there exists a deterministic continuous function
$\bar v^{ij,0}$ of class $\Pi_g$ such that for any $t\leq T$ and
$s\in [t,T]$, $\bar Y^{ij,0}_s=\bar v^{ij,0}(s,\x_s)$. Then, taking
$s=t$ and the limit as $n \rw \infty$ to obtain
$$
\underline{v}^{ij}(t,x):=\lim_{n\rw
\infty}\underline{v}^{ij,n}(t,x)\le \lim_{n\rw \infty}
v^{ij,n,0}(t,x)=\bar v^{ij,0}(t,x),\,\,\forall (t,x)\in \esp.
$$
But, $\bar v^{ij,0}$ and $\underline{v}^{ij,n}$ belong to $\Pi_g$
and $\underline{v}^{ij,n}\leq \underline{v}^{ij,n+1}$. Thus,
$\underline{v}^{ij}$ belongs to $\Pi_g$, for any $(i,j)\in \gam$.
In the same way one can show that $\bar v^{ij}$ belongs to $\Pi_g$,
for any $(i,j)\in \gam$. The last inequality follows from
(\ref{ineqvs}) and the definitions of $\bar v^{ij}$ and
$\underline{v}^{ij}$. \qed

\section{A comparison result}
In this section we investigate some qualitative properties of viscosity solutions of the system (\ref{main-systvi}). In particular, we show
in Corollary \ref{unicite} below that if the system (\ref{main-systvi}) admits a viscosity solution in the class $\Pi_g$, then it is unique and continuous.

\medskip\noindent
We first show that if $(u^{ij})_{\ij}$ (resp. $(w^{ij})_{\ij}$) is
an $usc$ subsolution (resp. $lsc$ supersolution) of
(\ref{main-systvi}) which belongs to $\Pi_g$, then for any $\ij$,
$u^{ij}\leq w^{ij}$. To begin with, we give an intermediary result
which is required in the second step of the proof of the comparison
result.

\begin{lem}\label{rslt_inter}
Let $\vec{u}:=(u^{ij})_{\ij}$ (resp. $\vec{w}:=(w^{ij})_{\ij}$) be
an $usc$ subsolution (resp. $lsc$ supersolution) of
(\ref{main-systvi}). Let $\tx$ and let $\tilde \Gamma
(t,x)$ be the following set.
$$
\tilde \Gamma (t,x):=\{\ij,
u^{ij}(t,x)-w^{i,j}(t,x)=\max_{(k,l)\in
\gam}\{u^{kl}(t,x)-w^{kl}(t,x)\}.
$$
Then, there exists $(i_0,j_0)\in\tilde \Gamma (t,x)$ such that
\be\lb{eqbord}
u^{i_0j_0}(t,x)>L^{i_0j_0}[\vec{u}](t,x)\,\,\, \mbox{ and
}\,\,\, w^{i_0j_0}(t,x)<U^{i_0j_0}[\vec{w}](t,x).\ee
\end{lem}

\noindent \pr: Let $\tx$ be fixed. Since the set $ \gam $ is finite
then $\tilde \Gamma (t,x)$ is not empty. Next, let us show the claim
by contradiction. So for any $(i,j) \in \tilde \Gamma (t,x)$ either
$u^{ij}(t,x)\leq L^{ij}[\vec{u}](t,x)$ or $w^{ij}(t,x)\geq
U^{ij}[\vec{w}](t,x)$ holds. Let us first assume that:
\begin{equation}\label{firstassump} \displaystyle{u^{ij}(t,x)\leq L^{ij}[\vec{u}](t,x).}
\end{equation}
 Then, there exists some $k\in \gami$, such that
$$
u^{ij}(t,x)\leq L^{ij}[\vec{u}](t,x)=u^{kj}(t,x)-\underline
g_{ik}(t,x).$$  But, since $\vec{w}$ is a supersolution to
(\ref{main-systvi}) we also deduce
$$
w^{ij}(t,x)\geq w^{kj}(t,x)-\underline g_{ik}(t,x),$$ which
implies that
$$
u^{ij}(t,x)-u^{kj}(t,x)\leq-\underline g_{ik}(t,x)\leq
w^{ij}(t,x)-w^{kj}(t,x)
$$ and then
$$
u^{ij}(t,x)-w^{ij}(t,x)\leq u^{kj}(t,x)-w^{kj}(t,x).
$$
Since, by assumption $(i,j)\in \tilde \Gamma(t,x)$, the two previous
inequalities are instead equalities, $(k,j)$ belongs to $\tilde
\Gamma(t,x)$, $k\neq i$ and finally it holds that
\be\label{equationuw1}
u^{ij}(t,x)-u^{kj}(t,x)=  -\underline g_{ik}(t,x)=w^{ij}(t,x)-w^{kj}(t,x).
\ee

\no Next, if (\ref{firstassump}) does not hold, then necessarily
$u^{ij}(t,x)> L^{ij}[\vec{u}](t,x)$ and $w^{ij}(t,x)\geq
U^{ij}[\vec{w}](t,x)$. As $(u^{ij})_{i,j}$ {is a subsolution of
(\ref{main-systvi}) we obtain
$$
u^{ij}(t,x) \le U^{ij}[\vec{u}](t,x)\leq u^{il}(t,x)+\bar g_{jl}(t,x), \quad
l\in {\Gamma_2}^{-j}.
$$
On the other hand, for some index $l\in {\Gamma_2}^{-j}$ it holds
$$
w^{ij}(t,x)- w^{il}(t,x)\geq \bar g_{jl}(t,x)\geq
u^{ij}(t,x)-u^{il}(t,x).
$$
Once more as $(i,j)\in \tilde \Gamma (t,x)$ then the
previous inequalities yield that $(i,l)\in \tilde \Gamma (t,x)$, $l\neq j$
and
\be\label{equationuw2}
u^{ij}(t,x)-u^{il}(t,x)=\bar g_{jl}(t,x)=w^{ij}(t,x)-
w^{il}(t,x). \ee
Repeating now this reasoning as many times as
necessary, and since $\gam$ is finite, there exits a loop
$(i_1,j_1),\ldots,(i_{N-1},j_{N-1}),(i_{N},j_{N})=(i_{1},j_{1})$ such
that
$$
\sum_{q=1,N-1}\varphi_{i_q,i_{q+1}}(t,x)=0,
$$
which contradicts Assumption (H4), whence the claim is proved. \qed \ms

\ms\no We now give the main result of this subsection.

\begin{thm}\label{comparison}
Assume that $\vec{u}$ = ($u^{ij}$)$_{(i,j) \in \Gamma^{1}\times
\Gamma^{2}}$ (resp. $\vec{w}$ = ($w^{ij}$)$_{(i,j) \in
\Gamma^{1}\times \Gamma^{2}}$) is an $usc$ (resp. $lsc$) subsolution
(resp. supersolution) of the system (\ref{main-systvi}) such that, for
any $\ij$, both $u^{ij}$ and $w^{ij}$ belong to $\Pi_{g}$, i.e.,
there exist two constants $\gamma$ and $C$ such that
 \begin{equation}\label{growthconditon}\forall \ij, \,\,
\forall (t,x)\in \esp, \,\,\, |u^{ij}(t,x)| +|w^{ij}(t,x)|\le
C(1+|x|^{\gamma}).
\end{equation}
Then, it holds that for any $\ij$,
$$
u^{ij}(t,x)\le
w^{ij}(t,x),\quad \tx.
$$
\end{thm}
\pr. Let us proceed by contradiction and let
($u^{ij})_{(i,j) \in \Gamma^{1}\times \Gamma^{2}}$ (resp. $\vec{w}$
= ($w^{ij}$)$_{(i,j) \in \Gamma^{1}\times \Gamma^{2}}$) be
$usc$ (resp. $lsc$) and a subsolution (resp. a supersolution) of the system
(\ref{main-systvi}) such that there exists $\eps_0>0$ and
$(t_0,x_0)\in \esp$ such that
\begin{equation}\label{contradiction} \max_{i,j} \big((u^{ij}-w^{ij})(t_0,x_0) \big) \geq
\epsilon_{0}.
\end{equation}
Next, w.l.o.g. we may assume that for any $\ij$,
\be\lb{inlim}
\lim_{|x|\rw \infty}(u^{ij}-w^{ij})(t,x)=-\infty. \ee
Indeed, if this is not the case, one may replace $w^{ij}$ with  $w^{ij, \vartheta,\mu}$ defined
by
$$
w^{ij, \vartheta,\mu}(t,x) = w^{ij}(t,x) + \vartheta e^{-\mu t}|x|^{2\gamma +2}, \tx,
$$
which is still an $usc$ supersolution of (\ref{main-systvi}) for
$\vartheta>0$ and $\mu\geq \mu_0$ which satisfies
(\ref{inlim}) (a proof of the supersolution property for good
choices of $\vartheta$ and $\mu$ can be found  in e.g. Pham (2009)
((\cite{Pham2009}), pp.76). Therefore, it suffices to show that
$u^{ij}(t,x)\leq w^{ij, \vartheta,\mu}(t,x),\,\,\, \tx,$ since, by
taking the limit as $\vartheta \rw 0$, we deduce that
$u^{ij}(t,x)\le w^{ij}(t,x),\,\,\, \tx.$ Thus, assume that
(\ref{contradiction}) and (\ref{inlim}) are satisfied. Then, there
exists $R>0$ such that
\begin{equation}\label{localizationmaxpoint}\ba{ll}
\max_{ (t,x) \in [0, T] \times \mathbb{R}^{k}}\max_{i,j} \{(u^{ij}
-w^{ij})(t,x)\}&=
 \max_{ (t,x) \in [0, T] \times B(0,R)}\max_{i,j} \{(u^{ij} -w^{ij})(t,x)\}\\{}&=\max_{ij}
  (u^{ij} -w^{ij})(t^{*}, x^{*})\geq \epsilon_{0} >0,
\ea
\end{equation}
where, $(t^{*}, x^{*})\in [0,T) \times B(0,R)$, where, $B(0,R)$ denotes the ball in $\R^ k$ with center the origin and radius $R$, since by definition
$u^{ij}(T,x)\leq {w}^{ij}(T,x)$, for all $\ij$. \ms

\medskip\noindent The remaining of the proof is obtained in two steps: the first step which is the main one establishes the comparison result under the additional condition (\ref{monotonicity-assumpt}) and the second step provides the proof in the general case. \ms

\noindent \paragraph*{\underline{Step 1}.} Let us make the following
assumption on the functions $(f^{ij})_{\ij}$. For all
$$
\lambda  >c(f^{i,j})(\Lambda-1),\; \ij,  \; (t,x,\vec{y},z)\in \esp\times \R^{\L+d},
\,\mbox{ and }\, (u,v)\in \R^2\,\mbox{s.t.}\, u \ge v,\\
$$
\begin{equation}\label{monotonicity-assumpt}\ba{l}
f^{ij}(t,x,[\vec{y}^{-(ij)},u],z)
-f^{ij}(t,x,[\vec{y}^{-(ij)},v],z) \le -\lambda(u-v),\ea
\end{equation}
where, $c(f^{ij})$ is the Lipschitz constant of
$f^{ij}$ w.r.t. $(y^{kl})_{(k,l)\in \gam}$. Next, let $(i_0,j_0)$ be an element of $\tilde \Gamma (t^*,x^*)$ that
satisfies (\ref{eqbord}). For $n\geq 1$, let $ \Phi^{i_{0},
j_{0}}_{n}$ be the function defined as follows.
$$
\Phi^{i_{0}, j_{0}}_{n}(t,x,y):=(u^{i_{0}j_{0}}(t,x)-w^{i_{0}
j_{0}}(t,y))-\phi_n(t,x,y),\quad (t,x,y)\in [0,T]\times \R^{k+k},
$$
where,
$$\phi_n(t,x,y):=n|x-y|^{2\g+2}+|x-x^*|^{2\g+2}+(t-t^*)^2.
$$
Since $\Phi^{i_{0}, j_{0}}_{n}$ is $usc$ in $(t,x,y)$, there
exists $(t_n,x_n,y_n)\in [0,T]\times B(0,R)^2$ such that
$$
\Phi^{i_{0}, j_{0}}_{n}(t_n,x_n,y_n)=\max_{(t,x,y)\in [0,T]\times
B(0,R)^2}\Phi^{i_{0}, j_{0}}_{n}(t,x,y).
$$
Moreover,
\be\begin{array}{l} \label{estiviscoge}\Phi^{i_{0},
j_{0}}_{n}(t^*,x^*,x^*)=u^{i_{0}, j_{0}}(t^*,x^*)-w^{i_{0},
j_{0}}(t^*,x^*) \leq u^{i_{0}, j_{0}}(t^*,x^*)-w^{i_{0},
j_{0}}(t^*,x^*)+ \phi_n(t_n,x_n,y_n) \\\qq \qq \qq \q \leq u^{i_{0},
j_{0}}(t_n,x_n)-w^{i_{0}, j_{0}}(t_n,y_n).\end{array}\ee
The definition of $\phi_n$ together with the growth condition of
$u^{ij}$ and $w^{ij}$ implies that $(x_n-y_n)_{n\geq 1}$ converges
to $0$. Next, for any subsequence $((t_{n_l},x_{n_l},y_{n_l}))_l$
which converges to $(\tilde t, \tilde x, \tilde x)$ we deduce from
(\ref{estiviscoge}) that
$$
u^{i_0j_0}(t^*,x^*)-w^{i_0j_0}(t^*,x^*)\leq u^{i_0,j_0}(\tilde
t,\tilde x)-w^{i_0j_0}(\tilde t,\tilde x),
$$ since $u^{i_0j_0}$ is $usc$ and $w^{i_0j_0}$ is $lsc$.
As the maximum of $u^{i_0j_0}-w^{i_0j_0}$ on $[0,T]\times B(0,R)$
is reached in $(t^*,x^*)$, then this last inequality is actually an
equality. Using the definition of $\phi_n$ and (\ref{estiviscoge}),
we deduce that the sequence $((t_{n},x_{n},y_{n}))_n$ converges to
$(t^*,x^*,x^*)$ which also implies
$$
n|x_n-y_n|^{2\g+2}\rightarrow 0 \,\,\,\mbox{ and
}\,\,\,(u^{i_0j_0}(t_n,x_n),w^{i_0j_0}(t_n,y_n))\rightarrow
(u^{i_0j_0}(t^*,x^*),w^{i_0j_0}(t^*,x^*)),
$$
as $n\to\infty$. This latter convergence holds, since we first
obtain from (\ref{estiviscoge}) and in taking into account that
$u^{i_0j_0}$ and $w^{i_0j_0}$ are $lsc$ and $usc$ respectively,
$$\ba{ll}
u^{i_0j_0}(t^*,x^*)-w^{i_0j_0}(t^*,x^*)&\leq \liminf_n (u^{i_0j_0}(t_n,x_n)-w^{i_0j_0}(t_n,y_n))\leq \limsup_n(u^{i_0j_0}(t_n,x_n)-w^{i_0j_0}(t_n,y_n))\\{}&\le
\limsup_nu^{i_0j_0}(t_n,x_n)-\liminf_nw^{i_0j_0}(t_n,y_n))\le u^{i_0j_0}(t^*,x^*)-w^{i_0j_0}(t^*,x^*).\ea
$$Thus the sequence $(u^{i_0j_0}(t_n,x_n)-w^{i_0j_0}(t_n,y_n))_{n\geq 0}$ converges
to $u^{i_0j_0}(t^*,x^*)-w^{i_0j_0}(t^*,x^*)$ as $n\rw \infty$. Then
$$
\ba{ll} \liminf_n
u^{i_0j_0}(t_n,x_n)=&u^{i_0j_0}(t^*,x^*)-w^{i_0j_0}(t^*,x^*)+\liminf_n
w^{i_0j_0}(t_n,y_n)\ge u^{i_0j_0}(t^*,x^*)\ge \limsup_n
u^{i_0j_0}t_n,x_n). \ea$$ It follows that the sequence
$(u^{i_0j_0}(t_n,x_n))_n$ converges to $u^{i_0j_0}(t^*,x^*)$ and
then $(w^{i_0j_0}(t_n,y_n))_n$ converges also to
$w^{i_0j_0}(t^*,x^*)$.\\ Next, recalling that $u^{i_{0} j_{0}}$
(resp. $w^{i_{0}j_{0}}$) is $usc$ (resp. $lsc$) and satisfies
(\ref{eqbord}), then, for $n$ large enough and at least for a
subsequence which we still index by $n$, we obtain
 \be \lb{ineqbord1} u^{i_{0} j_{0}}(t_n,x_n)>\max_{k\in
\gami}(u^{kj_0}(t_n,x_n)-g_{i_0k}(t_n,x_n)),\ee
 and \be
\lb{ineqbord2}w^{i_{0}j_{0}}(t_n,x_n)<\min_{l\in
\gamj}(w^{i_0l}(t_n,x_n)-g_{j_0l}(t_n,x_n)).
\ee
Applying now Crandall-Ishii-Lions's Lemma (see e.g. \cite{Crandalllions92} or \cite{Flemingetsoner}, pp.216) with
$\Phi^{i_{0}j_{0}}_{n}$ and $\phi_n$ at the point $(t_n,x_n,y_n)$
(for $n$ large enough in such a way that this latter triple will
belong to $[0,T]\times B(0,R)^2$), there exist $(p^n_u,q^n_u,M^n_u)\in \bar
J^{2,+}(u^{i_0j_0})(t_n,x_n)$ and $( p^n_w,q^n_w,M^n_w)\in \bar
J^{2,-}(w^{i_0j_0})(t_n,y_n)$ such that

$p^n_u-p^n_w=\partial_t\tilde \varphi_n(t_n,x_n,y_n)=2(t_n-t^*),\,$ $q^n_u \,\,(\mbox{resp. }q^n_w)\,=\partial_x\varphi_n(t_n,x_n,y_n)$
$(\mbox{resp.}-\partial_y
\varphi_n(t_n,x_n,y_n))$  and
\be \label{cil31}
\left (\begin{array}{ll} M_u^n&0\\
0&-N_w^n\end{array}\right )\leq A_n+\frac{1}{2n}A_n^2,
\ee
where, $A_n=D^2_{(x,y)}\varphi_n(t_n,x_n,y_n)$. Taking into account
that $(u^{ij})_{i,j}$ (resp. $(w^{ij})_{i,j}$) is a subsolution
(resp. supersolution) of (\ref{main-systvi}) and using once more
(\ref{ineqbord1}) and (\ref{ineqbord2}) we get
$$
-p^n_u-b(t_n,x_n)^\top.q^n_u-\frac{1}{2}
Tr[(\sigma\sigma^\top)(t_n,x_n)M^n_u]-f^{i_0j_0}(t_n,x_n,(u^{ij}(t_n,x_n))_{\ij},\sigma
(t_n,x_n)^\top.q^n_u)\leq 0,
$$
and
$$-p^n_w-b(t_n,y_n)^\top.q^n_w-\frac{1}{2}
Tr[(\sigma\sigma^\top)(t_n,y_n)M^n_w]-f^{i_0j_0}(t_n,y_n,(w^{ij}(t_n,y_n))_{\ij},\sigma
(t_n,y_n)^\top.q^n_w)\geq 0.
$$
Taking the difference between these two inequalities yields
\be\label{cil42}\begin{array}{l}
-(p^n_u-p^n_w)-(b(t_n,x_n)^\top .q^n_u-b(t_n,y_n)^\top . q^n_w)-
\frac{1}{2}
Tr[\{\sigma\sigma^\top(t_n,x_n)M^n_u-\sigma\sigma^\top(t_n,y_n)M^n_w\}]\\\qquad
-\{f^{i_0j_0}(t_n,x_n,(u^{ij}(t_n,x_n))_{\ij},\sigma
(t_n,x_n)^\top\cdot q^n_u)\\\qq\qq\qq\qq\qq\qq
-f^{i_0j_0}(t_n,y_n,(w^{ij}(t_n,y_n))_{\ij},\sigma
(t_n,y_n)^\top . q^n_w)\}\leq 0, \end{array} \ee
and then
$$
\begin{array}{l}
-\{f^{i_0j_0}(t_n,x_n,(u^{ij}(t_n,x_n))_{\ij},\sigma
(t_n,x_n)^\top.q^n_u)\\\qq\qq\qq\qq\qq\qq
-f^{i_0j_0}(t_n,y_n,(w^{ij}(t_n,y_n))_{\ij},\sigma
(t_n,y_n)^\top.q^n_w)\}\leq  \varrho_n, \end{array}
$$
with $\limsup_{n\rw \infty}\varrho_n\le 0$, using the fact that all
the terms in the first line of (\ref{cil42}) are converging
sequences. Linearizing $f^{i_0j_0}$ (see Appendix A3), which is
Lipschitz w.r.t. $(y^{ij})_{\ij}$, and using Assumption
(\ref{monotonicity-assumpt}), we obtain
$$
\l (u^{i_0j_0}(t_n,x_n)-w^{i_0j_0}(t_n,y_n))-\sum_{(i,j)\in \gam,
(i,j)\neq
(i_0j_0)}\Theta_n^{i,j}(u^{ij}(t_n,x_n)-w^{ij}(t_n,y_n))\leq
\varrho_n,
$$
where, $\Theta_n^{i,j}$ is the increment rate of $f^{i_0j_0}$ w.r.t.
$y^{ij}$ which is uniformly bounded (w.r.t. $n$) and non-negative thanks to the monotonicity assumption (H2). Therefore,
$$
\ba{ll}
\l (u^{i_0j_0}(t_n,x_n)&-w^{i_0j_0}(t_n,y_n))\\
{}&\leq \sum_{(i,j)\in \gam, (i,j)\neq
(i_0,j_0)}\Theta_n^{i,j}(u^{ij}(t_n,x_n)-w^{ij}(t_n,y_n))+
\varrho_n\\
{}&\leq c(f^{i_0j_0})\times \sum_{(i,j)\in \gam, (i,j)\neq
(i_0,j_0)}((u^{ij}(t_n,x_n)-w^{ij}(t_n,y_n))^++ \varrho_n. \ea
$$
Taking the limit as $n\rw \infty$ we obtain
$$
\ba{ll}
\l (u^{i_0j_0}(t^*,x^*)&-w^{i_0j_0}(t^*,y^*))\\
{}&\leq \limsup_{n\rw \infty}c(f^{i_0j_0})\left(\sum_{(i,j)\in
\gam, (i,j)\neq
(i_0,j_0)}(u^{ij}(t_n,x_n)-w^{ij}(t_n,y_n))^+\right)\\

{}& \le c(f^{i_0j_0})\left( \sum_{(i,j)\in \gam, (i,j)\neq
(i_0,j_0)}\limsup_{n\rw
\infty}(u^{ij}(t_n,x_n)-w^{ij}(t_n,y_n))^+\right)\\
{}& \le c(f^{i_0j_0})\left( \sum_{(i,j)\in \gam, (i,j)\neq
(i_0,j_0)}(u^{ij}(t^*,x^*)-w^{ij}(t^*,x^*))^+ \right)\ea,
$$
since $u^{ij}$ (resp. $w^{ij}$) is $usc$ (resp. $lsc$). As
$(i_0,j_0)$ belongs to $\tilde \Gamma (t^*,x^*)$, we obtain
$$
\l (u^{i_0j_0}(t^*,x^*)-w^{i_0j_0}(t^*,y^*))\leq c(f^{i_0j_0})
 \left( (\Lambda -1)(u^{i_0j_0}(t^*,x^*)-w^{i_0j_0}(t^*,x^*))\right).
 $$
But this is contradictory with 
(\ref{localizationmaxpoint}) and (\ref{monotonicity-assumpt}). Thus, for any $\ij$, $u^{ij}\leq
w^{ij}$. $\Box$ \ms

\noindent \paragraph*{\underline{Step 2}.} The general case.

\noindent For arbitrary $\l\in \R$, if $(u^{ij})_{\ij \in \gam}$ (resp.
$(w^{ij})_{\ij}$) is a subsolution (resp. supersolution)
of (\ref{main-systvi}) then $\tilde u^{ij}(t,x)=e^{\l
t}u^{ij}(t,x)$ and $\tilde w^{ij}(t,x)=e^{\l t}w^{ij}(t,x)$ is a
subsolution (resp. supersolution) of the following system of
variational inequalities with oblique reflection. For every $\ij$,
\be \label{sysvi3} \left\{
\begin{array}{l}
\min \{\tilde v^{ij}(t,x)-\max_{l\in \gami}\{\tilde
v^{lj}(t,x)-e^{\l t}\underline{g}_{il}(t,x)\};\\
\,\,\qq  \max\{\tilde v^{ij}(t,x)-\min_{k \in \gamj }\{e^{\l t}\bar
g_{jk}(t,x) + \tilde v^{ij}(t,x)\};-\partial_t\tilde v^{ij}(t,x)-
{\cal L}\tilde v^{ij}(t,x)+\l \tilde v^{ij}(t,x)\\\qq \qq \qq \qq
-e^{\l t}f^{ij}(t,x,(e^{-\l t} \tilde v^{ij}(t,x))_{\ij},e^{-\l
t}\sigma^\top(t,x).D_x \tilde v^{ij}(t,x))\}=0,\\ \tilde
v^{ij}(T,x)=e^{\l T}h^{ij}(x).
\end{array}\right.
\ee
But, by choosing $\l$ large enough the functions
$$
F^{ij}(t,x,(u^{kl})_{(k,l)\in \gam},z)=-\l u^{ij}+e^{\l t}f^{ij}(t,x,(e^{-\l
t}u^{kl})_{(k,l)\in \gam},e^{-\l t}z),\,\,\,\ij,
$$
satisfy Condition (\ref{monotonicity-assumpt}). Hence, thanks to the result stated in
Step 1, we have $\tilde u^{ij}\le \tilde v^{ij}$, $\ij$. Thus, $u^{ij}\le
v^{ij}$, for any $\ij$, which is the desired result.  \qed

\medskip\noindent As a consequence of this comparison result, we obtain the following one related to uniqueness
of the solution of (\ref{main-systvi}).
\begin{corollary}\label{unicite}
If the system (\ref{main-systvi}) admits a viscosity solution in the class $\Pi_g$, then it is unique and continuous.
\end{corollary}
\pr. Indeed, assume that $({v}^{ij})_{\ij}$ is a solution of
(\ref{main-systvi}) that belongs to $\Pi_g$. Then, thanks to
the previous comparison result, for any $\ij$ we have $v^{ij,*}\le
v^{ij}_*$. Thus, $v^{ij,*}=v^{ij}_*$ and then
$v^{ij}=v^{ij}_*=v^{ij,*}$, which means that $v^{ij}$ is
continuous. Next, if $(\hat v^{ij})_{\ij}$ is another
solution of (\ref{main-systvi}) in the class $\Pi_g$, then it is also
continuous and by the Comparison Theorem \ref{comparison} we have ${v}^{ij}\leq
\hat v^{ij}$  and ${v}^{ij}\geq \hat v^{ij}$. Hence,
$v^{ij}=\hat v^{ij},\,\,\,\ij $, i.e., uniqueness of
the solution of (\ref{main-systvi}). \qed

\section{Viscosity solution of the system (\ref{main-systvi})}
In this section we prove that the family
$(\bar{v}^{ij})_{i,j}$ constructed in Section 3 provides the unique continuous solution  in viscosity sense of the system
(\ref{main-systvi}). For sake of clarity, the proof is divided into several steps.

\begin{prop}\label{main-subsol}
The family $(\bar{v}^{ij})_{\ij}$ is a viscosity subsolution of the
system (\ref{main-systvi}).
\end{prop}
\pr. First recall that for each $\ij$, $\bar{v}^{ij}$ is
\textit{usc}, since $\bar v^{ij} = \lim _{m}\searrow v^{ij,m}$,
where $v^{ij,m}$ is a continuous function solution of the system
 (\ref{newauxiliary-systvi}). Thus, for any $\ij$, it holds that $\bar
 v^{ij,*}=\bar v^{ij}$. Next, at $T$ we have
$$
\bar v^{ij}(T,x) = \lim_m \searrow v^{ij,m}(T,x)=h^{ij}(x),\quad x\in \R^k.
$$
We shall now prove that, for any $(t,x)$ in $[0,T) \times
\R^k$, any $(\ij$ and ($\underline{p}, \underline{q}
,\underline{M}$) in $\bar{J}^{+}\bar{v}^{ij}(t,x)$,

\be  \label{subsolutionprop} \begin{array}{l} \min
[\big(\bar{v}^{ij}- \bar {L}^{ij}\big)(t,x),\max\{\big(\bar{v}^{ij}-
\bar{U}^{ij}\big)(t,x),  \qquad  \\\qquad-\underline{p} -b(t, x)
\underline{q}-\frac{1}{2}\textrm{Tr}(\sigma\sigma^{T}(t,x)
\underline{M}) -f^{ij}(t,x, (\bar v^{kl}(t,x))_{(k,l)\in \gam},
\sigma^\top (t,x).\underline q )\}] \le 0,\end{array}
 \ee
with $ \bar{L}^{ij} $ and $ \bar{U}^{ij}$ defined as follows:
$$
\displaystyle{ \bar{L}^{ij}(t,x) = \max_{k \in
(\Gamma^{1})^{-i}}\big( \bar{v}^{kj}(t,x)
-\underline{g}_{ik}(t,x)}\big) \;\,\, \mbox{and}\,\,
 \; \displaystyle{ \bar{U}^{ij}(t,x) =
 \min_{l \in (\Gamma^{2})^{-j}} \big(\bar{v}^{il}(t,x) +\bar{g}_{jl}(t,x)}\big).
 $$
Now, let $\ij$ be fixed. Then it is equivalent to show that,
either
  \begin{equation}\label{firstproperty}\big(\bar{v}^{ij} -\bar{L}^{ij}\big)(t,x) \le 0,  \end{equation}
or
 \begin{equation}\label{secondproperty}
\max\{ \big(\bar{v}^{ij} -\bar{U}^{ij}\big)(t,x), - \underline{p}-
b(t,x)\underline{q} -\frac{1}{2}\textrm{Tr}(\sigma \sigma^{T}(t,x)
\underline{M}) -f^{ij}(t,x, \vec{v}(t,x),\sigma^\top
(t,x).\underline q) \}
 \le 0.
   \end{equation}
  If (\ref{firstproperty}) is satisfied then the subsolution property (\ref{subsolutionprop})
  holds. Therefore, from now on, we suppose that there exists $\epsilon_{0} >0$ such that
  \begin{equation}\label{strictsubsolconstraint} \bar{v}^{ij}(t,x) \ge \bar{L}^{ij}(t,x) + \epsilon_{0},
  \end{equation}
 and show (\ref{secondproperty}). Thanks to the decreasing
 convergence of $(v^{ij,m})_{m\geq 0}$ to $\bar v^{ij}$, $\ij$,
 there exists $m_0$ such that for any $m\ge m_0$, we have
 \begin{equation}\label{auxiliaryconstraint}
 v^{ij,m }(t,x) \ge L^{ij,m}(t,x) + \frac{\epsilon_0}{2}. \end{equation}
Next, by continuity of $ v^{ij,m}$ and $ L^{ij,m}$, there exists a neighborhood $\mathcal{O}_{m}$ of $(t,x)$ such that
\begin{equation}\label{strictsubsolconstraint2} {v}^{ij,m}(t',x') \ge {L}^{ij,m}(t',x') +
\frac{\epsilon_{0}}{4},\quad (t',x')\in \cO_m.
  \end{equation}
Now, by  Lemma 6.1 in \cite{Crandalllions92} there
exists a subsequence $((t_{k},x_{k}))_{k\geq 0}$ such that
\begin{equation}\label{eq:unifconvergence}  (t_{k},x_{k})
\rightarrow_{k\to\infty} (t,x) \; \textrm{ and }\;\;\bar{v}^{ij}(t, x) =
\lim_{k\to \infty} v^{ij, k}(t_{k},x_{k}). \end{equation}
Moreover, there exists a sequence $(p_{k},q_{k},M_{k}) \in \bar{J}^{+}\big(v^{ij,
k}(t_{k},x_{k})\big)$ such that
\begin{equation}\label{superjetconvergence}
\lim_{k\to\infty}(p_{k},q_{k},M_{k})=(\underline{p}, \underline{q}, \underline{M}).
\end{equation}
But, the subsequence $((t_{k},x_{k}))_{k\geq 0}$ can be chosen in such a way
that for any $k\ge 0$, $(t_{k},x_{k})\in \cO_k$. Applying now the
viscosity subsolution property of $v^{ij, k}$ (which satisfies
(\ref{newauxiliary-systvi})) at $(t_{k},x_{k})$ and taking into
account of (\ref{strictsubsolconstraint2}) we obtain
\begin{equation}\label{eq:viscosityprop}
 -p_{k}-  b(t_{k},x_{k})^\top .q_k -\frac{1}{2}\textrm{Tr}(\sigma \sigma^\top(t_{k},x_{k}) M_{k}) -
 \bar f^{ij,k}(t_{k},x_{k}, (\bar v^{pq,k}(t_k,x_k))_{(p,q)\in \gam}, \sigma (t_{k},x_{k})^\top q_k )\leq 0,
 \end{equation}
where, once more,
$$
\bar{f}^{ij,k}(s, x, (y^{pq})_{(p,q)\in \gam},z) = f^{ij}
(s,x,(y^{pq})_{(p,q)\in \gam},z) -k\big(y^{ij} -\min_{l\in \gamj}
(y^{il} + \bar{g}_{jl}(s,x) )\big)^{+}.
$$
Next, thanks to the boundedness of the sequence $((t_{k},x_k))_{k\geq 0}$, the uniform polynomial
growth of $\bar v^{pq,k}\,\, k\geq 0$, (by Proposition
\ref{firstresult} and Corollary \ref{cndpol}), the assumptions
(H0)-(H2) on $b$, $\sigma$ and $f^{ij}$, and the convergence of
$((p_k,q_k,M_k))_k$ to $(\underline{p}, \underline{q},
\underline{M})$,  we deduce from (\ref{eq:viscosityprop}) that
$$ \epsilon_k:=\big(\bar v^{ij,k}(t_{k},x_{k}) -\min_{l \in
(\Gamma^{2})^{-j}}\big(\bar
v^{il,k}(t_{k},x_{k})+\bar{g}_{jl}(t_{k},x_{k})\big))^{+}
\rightarrow 0,\quad k\to\infty.
$$
But, for any fixed $(t,x)$ and $k_0$, the sequence $(\bar v^{il,k}(t,x))_{k \ge k_{0}}$ is
decreasing and then for any $k\ge k_{0}\ge m_0$,
$$
\displaystyle{ \bar v^{ij,k}(t_{k},x_{k}) \le \min_{l\neq j}\big(\bar v^{il,k}(t_{k},x_{k})
+\bar{g}_{jl}(t_{k},x_{k})\big) + \epsilon_{k} \le \min_{l\neq
j}\big(\bar v^{il,k_{0}}(t_{k},x_{k})
+\bar{g}_{jl}(t_{k},x_{k})\big) + \epsilon_{k}.}
$$
Taking now the limit as $k\rw +\infty$, in view of the continuity of $\bar v^{il,k_{0}}$, we get
$$
\lim_{k}\bar v^{ij, k}(t_{k},x_{k}) = \bar{v}^{ij}(t,x) \le
\min_{j\neq l} (v^{il, k_{0}}(t, x) +\bar{g}_{jl}(t, x)).
$$
Finally, passing to the limit as $k_{0}$ goes to $+\infty$ to obtain
\be \label{firstineq} \bar{v}^{ij}(t,x) \le \min_{j\neq l}(\bar{v}^{il}(t,x)
 +\bar{g}_{jl}(t, x))= \bar{U}^{ij}(t,x).
 \ee
Let us now consider a subsequence of $(k)$, which we denote by
$(k_l)$, such that for any $(p,q)\in \gam$, the sequence $(\bar
v^{pq,k}(t_{k_l},x_{k_l}))_l$ is convergent. This subsequence
exists since the functions $\bar v^{ij,{k_l}}$ are uniformly of
polynomial growth (by Proposition \ref{firstresult} and Corollary
\ref{cndpol}). Then, taking the limit w.r.t. $l$ in equation
(\ref{eq:viscosityprop}), we obtain
  \begin{equation}\label{eq:limit-prop}\ba{ll}
   - \underline{p}- \underline{q} b(t,x) -\frac{1}{2}\textrm{Tr}(\sigma \sigma^{\top}(t,x)
  \underline{M})& \le \limsup_{l\rw \infty}
 \bar f^{ij,k_l}(t_{k_l},x_{k}, (\bar v^{pq,k_l}(t_{k_l},x_{k_l}))_{(p,q)\in \gam}, \sigma (t_{k_l},x_{k_l})^\top
 q_{k_l})\\{}&\le \limsup_{l\rw \infty}
f^{ij}(t_{k_l},x_{k}, (\bar v^{pq,k_l}(t_{k_l},x_{k_l}))_{(p,q)\in
\gam}, \sigma (t_{k_l},x_{k_l})^\top
 q_{k_l})
\\{}&=
f^{ij}(t,x, (\lim_{l\rw \infty}\bar
v^{pq,k_l}(t_{k_l},x_{k_l}))_{(p,q)\in \gam}, \sigma (t,x)^\top
\underline q),
 \ea
\end{equation}
since $f^{ij}$ is continuous in $(t,x,\vec{y},z)$. Now for any
$(p,q)\in \gam$, since $\bar v^{pq,n}$ is continuous and $\bar
v^{pq,n}\geq \bar v^{pq,n+1}$, $\forall n\ge 0$, it holds that
$$
\bar v^{pq,*}(t,x)=\bar v^{pq}(t,x)=\limsup_{t'\rw t, x'\rw x,
n\rw \infty}\bar v^{pq,n}(t',x'), \qq  (t,x)\in [0,T]\times
\R^k.
$$
Therefore, for any $(p,q)\in \gam$ $((p,q)\neq (i,j))$, we have
\begin{equation}\label{eq:usc-prop}
\bar v^{pq}(t,x)\ge \lim_{l\rw \infty}\bar
v^{pq,k_l}(t_{k_l},x_{k_l})\mbox{ and }\bar v^{ij}(t,x)= \lim_{l\rw \infty}\bar
v^{ij,k_l}(t_{k_l},x_{k_l}).
\end{equation}
 As $f^{ij}$ is non-decreasing
w.r.t. $y^{kl}$, $(k,l)\in \gam, (k,l)\neq (i,j)$, we deduce from (\ref{eq:limit-prop}) and (\ref{eq:usc-prop})
that
\begin{equation}\label{eq:limit-prop2}\ba{l}
   - \underline{p}- \underline{q} b(t,x) -\frac{1}{2}\textrm{Tr}(\sigma \sigma^{\top}(t,x)
  \underline{M}) \le
f^{ij}(t,x, (\bar v^{pq}(t,x))_{(p,q)\in \gam}, \sigma (t,x)^\top
 \underline q).
 \ea
\end{equation}
Finally, under the condition (\ref{strictsubsolconstraint}), the
relations (\ref{eq:limit-prop2}), (\ref{firstineq}) imply that
(\ref{secondproperty}) is satisfied. Thus $\bar v^{ij}$ is a
viscosity subsolution for
$$
\left\{
\begin{array}{l} \min [\big(\bar{v}^{ij}-
\underline{L}^{ij}\big)(t,x),\max\{\big(\bar{v}^{ij}-
\underline{U}^{ij}\big)(t,x),  \qquad  \\\qquad-\underline{p} -b(t,
x)^\top \underline{q}-\frac{1}{2}\textrm{Tr}(\sigma\sigma^\top(t,x)
\underline{M}) -f^{ij}(t,x, (\bar v^{kl}(t,x))_{(k,l)\in \gam},
\sigma^\top (t,x).\underline q )\}] = 0,\\
\bar v^{ij}(T,x)= h^{ij}(x).
\end{array}\right.
 $$
Since $(i,j)$ is arbitrary, $(\bar v^{ij})_{i,j\in \gam}$ is a
viscosity subsolution for (\ref{main-systvi}). This finishes the proof.\qed

\begin{prop}\label{main-supersol} Let $m_0$ be fixed in $\mathbb{ N}$. Then, the family $(\bar{v}^{ij,m_0})_{\ij}$ is a viscosity supersolution
of the system (\ref{main-systvi}).
\end{prop}
\pr. First and thanks to Proposition \ref{firstresult}, we know that the
triples $((\bar Y^{ij,m_0},\bar Z^{ij,m_0},\bar
K^{ij,m_0}))_{(i,j)\in \gam}$ introduced in (\ref{multidimsystem}),
solve the following system of reflected BSDEs: For every $\ij$,
\be\label{mltdim}\left\{\begin{array}{l} \bar
Y^{ij,m_0}\in \cS^{2,1}, \,\,\bar Z^{ij,m_0} \in \cH^{2,d} \mbox{
and
}\bar K^{ij,m_0 } \in \cA^{2,+}\,\,;\\
\bar Y^{ij,m_0}_s=h^{ij}(X^{t,x}_T)+\int_{s}^{T}\bar
{f}^{ij,m_0}(r,X^{t,x}_r,
(\bar Y^{kl,m_0}_r)_{(k,l)\in \gam},\bar Z^{ij,m_0}_r)dr+\bar K^{ij,m_0}_T-\bar K^{ij,m_0}_s-
\int_s^T\bar Z^{ij,m_0}_rdB_r\\
\displaystyle{\bar Y^{ij,m_0}_{s}\geq \max_{k\in
(\Gamma^{1})^{-i}}\{\bar
Y^{kj,m_0}_s-\underline{g}_{ik}(s,X^{t,x}_s)\}},\,\,
s\leq T\,;\\
\int_0^T(\bar Y^{ij,m_0}_s-\max_{k\in (\Gamma^{1})^{-i}}\{\bar
Y^{kj,m_0}_s-\underline{g}_{ik}(s,X^{t,x}_s)\})d\bar
K^{ij,m_0}_s=0\end{array}\right.\ee
where, for any $(i,j) \in \gam$
and $(s, \vec{y}, z^{ij})$,
$$
\bar{f}^{ij,m_0}(s, X^{t,x}_s,\vec{y},z^{ij}) = f^{i,j}(s,X^{t,x}_s,(y^{kl})_{(k,l)\in \gam},z^{ij})
-m_0\big(y^{ij} -\min_{l\in \gamj} (y^{il} + \bar{g}_{jl}
(s,X^{t,x}_s))\big)^{+}.
$$
Furthermore, it holds true that
$$
\forall \ij,\,\,\, \forall \tx, \,\,\,\forall s\in
[t,T],\,\,\,\bar Y^{ij,m_0}_s=\bar v^{ij,m_0}(s,X_s^{t,x}).
$$
On the other hand, we note that $\bar Y^{ij,m_0}$ is the value
function of a zero-sum Dynkin game (see Appendix A4), i.e., it
satisfies, for all $s\leq T$ \be\label{eqconv2}\ba{l} \bar
Y^{ij;m_0}_s=\esssup_{\sigma \geq s}\essinf_{\tau \geq
s}\E[\int_s^{\sigma \wedge \tau}f^{ij}(r,X^{t,x}_r,(\bar
Y^{ij;m_0}_r)_{\ij},\bar Z^{ij;m_0}_r)dr+\\\qq \max_{k\in
\gami}\{\bar Y^{kj;m_0}_\sigma-\underline{g}_{ik}(\sigma,
X^{t,x}_\sigma)\}\ind_{[\sigma <\tau]}+\{\bar Y^{ij;m_0}_\t\wedge
\min_{l\in \gamj}\{\bar Y^{il;m_0}_\tau-\bar{g}_{jl}(\tau,
X^{t,x}_\tau)\}\}\ind_{[\tau\leq\sigma <T]}\\\qq\qq\qq\qq
+h^{ij}(X^{t,x}_T)\ind_{[\tau=\sigma =T]}|\cF_s].\ea \ee Thus, by
Theorem 3.7 in Hamad\`ene-Hassani (2005) (\cite{Hamhassanni}), it
follows that $\bar v^{ij,m_0}$ is the unique viscosity solution for
the following PDE with two obstacles.
$$\left\{
\begin{array}{l} \min [\vartheta(t,x)-\max_{k\in
\gami}\{\bar v^{kj,m_0}(t,x)-\underline g_{ik}(t,x)\},
\max\{\vartheta(t,x)-\bar v^{ij,m_0}(t,x)\vee \min_{l\in
\gamj}(\bar v^{il,m_0}(t,x)-\bar g_{jl}(t,x)), \\\qq\q -\partial_t
\vartheta(t,x)-b(t, x)^\top
D_x\vartheta(t,x)-\frac{1}{2}\textrm{Tr}(\sigma\sigma^\top (t,x)
D^2_{xx}\vartheta(t,x)) -\\ \qq\qq\q f^{ij}(t,x, (\bar
v^{kl,m_0}(t,x))_{(k,l)\in \gam},
\sigma^\top (t,x).D_x\vartheta(t,x))\}] = 0\,;\\
\vartheta(T,x)= h^{ij}(x).
\end{array}\right.
 $$
Next, let $\tx$ and $(p,q,M)\in {\bar J}^-[\bar
v^{ij,m_0}](t,x)$. As $\bar v^{ij,m_0}$ is a solution in a
viscosity sense of the previous PDE with two obstacles then it
holds that
\be\lb{ineq1}\bar v^{ij,m_0}(t,x)\geq \max_{k\in
\gami}\{\bar v^{kj,m_0}(t,x)-\underline g_{ik}(t,x)\}
\ee
and
\be\lb{ineq2}\ba{l} \max\{\bar v^{ij,m_0}(t,x)-\bar
v^{ij,m_0}(t,x)\vee \min_{l\in \gamj}(\bar v^{il,m_0}(t,x)-\bar
g_{jl}(t,x))\,;\\\qq \qq -p-b(t, x)^\top
q-\frac{1}{2}\textrm{Tr}(\sigma\sigma^\top (t,x) M) -f^{ij}(t,x,
(\bar v^{kl,m_0}(t,x))_{(k,l)\in \gam}, \sigma^\top (t,x).q)\}\geq
0.\ea
\ee
But, for any constants $a,b\in \R$ we have $a -(a\vee b) \le
a-b$ and thus $ a-(a\vee b) \ge 0 \Rightarrow a-b \ge 0.$ Therefore,
we have from (\ref{ineq2}), $$\ba{l} \max\{\bar v^{ij,m_0}(t,x)-
\min_{l\in \gamj}(\bar v^{il,m_0}(t,x)-\bar g_{jl}(t,x))\,;\\\qq
\qq -p-b(t, x)^\top q-\frac{1}{2}\textrm{Tr}(\sigma\sigma^\top (t,x)
M) -f^{ij}(t,x, (\bar v^{kl,m_0}(t,x))_{(k,l)\in \gam},
\sigma^\top (t,x).q)\}\geq 0.\ea
$$
Combining this inequality with
(\ref{ineq1}) and since $v^{ij,m_0}(T,x)= h^{ij}(x)$  it follows
that $v^{ij,m_0}$ is a viscosity supersolution of the system
$$
\left\{
\begin{array}{l} \min [\vartheta(t,x)-
\max_{k\in \gami}(\bar v^{kj,m_0}(t,x)-\underline g_{ik}(t,x))
\,;\,\max\{\vartheta(t,x)- \min_{l\in \gamj}(\bar
v^{il,m_0}(t,x)-\bar g_{jl}(t,x)); \qquad
\\\qquad-\partial_t \vartheta(t,x)-b(t, x)^\top
D_x\vartheta(t,x)-\frac{1}{2}\textrm{Tr}(\sigma\sigma^\top (t,x)
D^2_{xx}\vartheta(t,x)) -\\\qq \qq f^{ij}(t,x, (\bar
v^{kl,m_0}(t,x))_{(k,l)\in \gam},
\sigma^\top (t,x).D_x\vartheta(t,x))\}] = 0,\\
\vartheta(T,x)= h^{ij}(x).
\end{array}\right.
$$
Since $(i,j)$ is arbitrary in ${\gam}$, the system of continuous
functions $(v^{ij,m_0})_{\ij}$ is a supersolution of
(\ref{main-systvi}). \qed

\ms\no Consider now the set $\cU_{m_0}$ defined as follows.
$$
\mathcal{U} = \{\vec{u}:=(u^{ij})_{\ij}\; \textrm{s.t.} \; \vec{u} \; \textrm{is a subsolution of
(\ref{main-systvi})} \; \textrm{and} \;\forall \,(i,j)\in \gam, \,\,
 \bar{v}^{i,j}\le u^{i,j} \le \bar v^{ij,m_{0}}\}.
 $$
Obviously, $\cU_{m_0}$ is not empty since it contains $(\bar
v^{ij})_{\ij}$. Next for $\tx$ and $\ij$, let us set:
$$
^{m_0}v^{ij}(t,x)=\sup\{u^{ij}(t,x), \,\,(u^{kl})_{(k,l)\in \gam}\in \cU_{m_0}\}.
$$
We now state the main result of this section.
\begin{thm}\label{main1} The family $(^{m_0}{v}^{ij})_{\ij}$ does not depend
on $m_0$ and is the unique continuous viscosity solution in the class $\Pi_g$ of the system
(\ref{main-systvi}).
\end{thm}
\pr. We first note that for any $\ij$, $\bar v^{ij}\le \,^{m_0}\!v^{i,j}\le \bar v^{ij,m_0}$. Since $\bar v^{ij}$ and $\bar
v^{ij,m_0}$ are of polynomial growth, then $(^{m_0}v^{ij})_{\ij}$ belongs to $\Pi_g$.

\ms\no The remaining of the proof is divided into two steps and mainly consists in adapting the Perron's method
 (see Crandall-Ishii-Lions, \cite{Crandalllions92} Theorem 4.1, pp 23) to construct a viscosity solution to our general system of PDEs. To ease notation, we denote $(^{m_0}{v}^{ij})_{\ij}$ by
$({v}^{ij})_{\ij}$ as no confusion is possible.

\paragraph{Step 1.} We first show that $(v^{ij})_{\ij}$ is a subsolution of (\ref{main-systvi}). Indeed, it is clear that for any $\tx$, $\bar v^{ij}(t,x)\leq v^{ij}(t,x)\le \bar v^{ij,m_0}(t,x)$. This implies that $\bar
v^{ij}\leq v^{ij,*}\le \bar v^{ij,m_0}$ since, as pointed out
previously, $\bar v^{ij}$ is $usc$ and $v^{ij,m_0}$ is continuous.
Therefore, for any $x\in \R^k$, we have $v^{ij,*}(T,x)=h^{ij}(x)$,
since $\bar v^{ij}(T,x)=\bar v^{ij,m_0}(T,x)=h^{ij}(x)$. \ms

\no Next, fix $\ij$ and let $(\tilde{v}^{ij})_{\ij}$ be an arbitrary element of
$\mathcal{U}_{m_0}$. Then, for any $(t,x)\in
[0,T)\times \R^k$ and any $(p,q,M)\in \bar J^+{\tilde
v}^{i,j,*}(t,x)$ we have
$$
\left\{
\begin{array}{l} \min [\big(\tilde{v}^{ij,*}-
L^{ij}((\tilde{v}^{kl,*})_{k,l})\big)(t,x),\max\{\big(\tilde{v}^{ij,*}-
U^{ij}((\tilde{v}^{kl,*})_{k,l}\big)(t,x)\,,  \qquad  \\\qquad -p
-b(t, x)^\top q-\frac{1}{2}\textrm{Tr}(\sigma\sigma^{\top}(t,x) M)
-f^{ij}(t,x, (\tilde{v}^{kl,*}(t,x))_{(k,l)\in \gam},
\sigma^\top (t,x)q )\}] \le 0.\\
\end{array}\right.
 $$
By definition we have $\tilde{v}^{kl} \le v^{kl}$ and then
$\tilde{v}^{kl,*} \le v^{kl,*} $ for any $(k,l)\in \gam$. Since the operators $ \vec{w}= (w^{kl})_{k,l}
\mapsto \tilde{v}^{ij*} -L^{i,j}((w^{kl})_{k,l})$,
$\vec{w}=(w^{kl})_{k,l} \mapsto \tilde{v}^{ij*}
-U^{ij}((w^{kl})_{k,l})$ are decreasing, in view of the monotonicity property (H2) of
the generator $f^{ij}$, we have
$$\left\{
\begin{array}{l} \min [\big(\tilde{v}^{ij,*}-
L^{ij}((v^{kl,*})_{k,l})\big)(t,x)\,;\q
\max\{\big(\tilde{v}^{ij,*}- U^{ij}((v^{kl,*})_{k,l})\big)(t,x),
\qquad  \\\qquad -p -b(t, x)^\top
q-\frac{1}{2}\textrm{Tr}(\sigma\sigma^\top(t,x) M)
-f^{ij}(t,x,[(v^{kl,*}(t,x))_{ (k,l)\in \gam \atop (k,l) \neq
(i,j)},\tilde v^{ij,*}(t,x)],
\sigma^\top (t,x)q )\}] \le 0\\
\end{array}\right.
 $$
where $[(v^{kl,*}(t,x))_{ (k,l)\in \gam \atop (k,l) \neq
(i,j)},\tilde v^{ij,*}(t,x)]$ is obtained from $(\tilde
v^{kl,*}(t,x))_{ (k,l)\in \gam})$ by replacing $\tilde
v^{ij,*}(t,x)$ with $v^{ij,*}(t,x)$. This means that $(t,x) \in
[0,T)\times \R^k \longrightarrow\tilde{v}^{i,j}(t,x )$ is a
subsolution of the following equation.
$$\left\{
\begin{array}{l} \min [\big( w-
L^{ij}((v^{kl,*})_{k,l})\big)(t,x),\q \max\{\big(w-
U^{ij}((v^{kl,*})_{k,l})\big)(t,x),  \qquad  \\\qquad
\displaystyle{-p -b(t, x)^\top
q-\frac{1}{2}\textrm{Tr}(\sigma\sigma^\top(t,x) M)
-f^{ij}(t,x,[(v^{kl,*}(t,x))_{ (k,l)\in \gam \atop (k,l) \neq
(i,j)},w], \sigma^\top (t,x)q )\}] =0}.
\end{array}\right.
 $$
Next, relying on the lower semi continuity of the function
$$
\left\{
\begin{array}{l}
(t,x,w,p, q, M)\in [0,T]\times \R^{k+1+1+k}\times \mathbb{S}^k \longmapsto
\q \min \{\left( w- \max_{k \neq i}(v^{kj,*}(t,x)  - \underline{g}_{ik}(t,x))\right),\\
\qq\qq \max [\left( w- \min_{l \neq j}(v^{il,*}(t,x) +
\bar{g}_{jl}(t,x) )\right), \displaystyle{-p -b(t,x)^\top
q-\frac{1}{2}\textrm{Tr}(\sigma\sigma^\top(t,x) M) }\\ \qq
\qq\qq\qq\qq\qq
 -f^{ij}(t,x,[(v^{kl,*}(t,x))_{ (k,l)\in \gam \atop (k,l) \neq (i,j)},w], \sigma^\top (t,x)q ) ]\}\\
\end{array} \right.
$$
and using Lemma 4.2, in Crandall {\it et al.} (1992) (\cite{Crandalllions92}, pp.23), related to suprema of subsolutions, combined
 with the above result, it holds that $v^{ij}$ is
a subsolution of the following equation:
\be\label{equationseparee}\left\{
\begin{array}{l} \min [\big( w-
L^{ij}((v^{kl,*})_{k,l})\big)(t,x),\q \max\{\big(w-
U^{ij}((v^{kl,*})_{k,l})\big)(t,x),  \qquad  \\\qquad
\displaystyle{-p -b(t, x)^\top
q-\frac{1}{2}\textrm{Tr}(\sigma\sigma^\top(t,x) M)
-f^{ij}(t,x,[(v^{kl,*}(t,x))_{ (k,l)\in \gam \atop (k,l) \neq
(i,j)},w(t,x)], \sigma^\top (t,x)q )\}] =0},\\
v^{ij}(T,x)=h^{ij}(x).
\end{array}\right.
 \ee
Since $(i,j)$ is arbitrary in $\gam$, $(v^{ij})_{\ij}$ is a
 subsolution of (\ref{main-systvi}).

\paragraph{Step 2.} In this step we use the so called Perron's method to show that $(v^{ij})_{\ij}$ is a viscosity
supersolution of (\ref{main-systvi}).

\ms\no Indeed, we first note for any $(i,j)\in \gam$,
$$\underline v^{ij}=\underline v^{ij}_*\leq \bar v^{ij}_*\leq
v^{ij}_* \le \bar v^{ij,m_0}_*=\bar v^{ij,m_0},
$$
since $\bar v^{ij,m_0}$ is continuous and $\underline v^{ij}$ is $lsc$.
Therefore, for any $x\in \R^k$,
\be \label{cdtterminale}
v^{ij}_*(T,x)=h^{ij}(x)\ee thanks to $\underline
v^{ij}(T,x)=h^{ij}(x)=\bar v^{ij,m_0}(T,x)$. Next, assume that $(v^{i,j})_{\ij}$ is not a supersolution for
(\ref{main-systvi}). Then, taking into account of
(\ref{cdtterminale}) and Remark \ref{allegement}, there exists at
least one pair ($i,j$) such that $v^{i,j}$ does not satisfy the viscosity supersolution property: this means that for some point $(t_0, \;x_0)\in
[0,T)\times \R^k$ there exists a triple ($p,q,M$) in
$\mathcal{{J}}^{-}(v^{ij}_{*})(t_0, x_0)$ such that
\begin{equation}\label{contradiction-2}
  \left\{
\begin{array}{l} \min [\big(v_{*}^{ij}-
L^{ij}((v_*^{kl})_{(k,l)\in
\gam})\big)(t_0,x_0)\,;\,\max\{\big(v_{*}^{ij}-
U^{ij}((v_*^{kl})_{(k,l)\in \gam})\big)(t_0,x_0), \qquad
\\\qquad -p -b(t_0, x_0)^\top q
-\frac{1}{2}\textrm{Tr}(\sigma\sigma^\top (t_0,x_0) M)
-f^{ij}(t_0,x_0, (v_{*}^{kl}(t_0,x_0))_{(k,l)\in \gam},
\sigma^\top (t_0,x_0).q)\}] < 0.
\end{array}\right.
\end{equation}
We now follow the same idea as in Crandall {\it et al.} (1992) (\cite{Crandalllions92}, pp.24). For any positive $\delta$, $\gamma$ and $ r$, set
$u_{\delta, \gamma}$ and $B_r$ as follows:
$$
\displaystyle{u_{\delta, \gamma}(t,x)  =
v^{ij}_{*}(t_0,x_0) +\delta +\langle{ q, x-x_{0} \rangle} +p(t-t_0)
+ \frac{1}{2}\langle{(M -2\gamma)(x-x_0), (x-x_0) \rangle}},
$$
and
$$
 B_{r} := \{\tx, \;\textrm{s.t.} \; |t-t_0| +|x-x_0| <r   \}.
$$
Using (\ref{contradiction-2}) and continuity of all the data, choosing $\d$, $\g$
small enough we obtain
\begin{equation}\label{upst0x0}
  \left\{
\begin{array}{l} \min [\big(v_{*}^{ij} +\delta -
L^{ij}((v_*^{kl})_{(k,l)\in
\gam}))(t_0,x_0),\,\,\max\{\big(v_{*}^{ij}+\delta-
U^{ij}((v_*^{kl})_{(k,l)\in \gam}))(t_0,x_0),   \\ \qquad  - p
-b(t, x)^\top q -\frac{1}{2}\textrm{Tr}(\sigma\sigma^\top (t_0,x_0)
(M -2\gamma)) \\\qquad \q -f^{ij}(t_0,x_0,
[(v_{*}^{kl}(t_0,x_0))_{(k,l) \neq
(i,j)},v_{*}^{ij}(t_0,x_0)+\delta],\sigma^\top (t_0,x_0)q)\}] < 0.\\
\end{array}\right.
\end{equation}
Next, let us define the function $\Upsilon$ as follows.
$$
\begin{array}{l}   \Upsilon(t,x) = \;
 \displaystyle{\min \left\{ u_{\delta, \gamma}(t,x) - \max_{k\neq i} \big(v_{*}^{kj} - \underline{g}_{ik}\big)(t,x), \;
  \max \{u_{\delta, \gamma}(t,x) - \min_{l \neq j} \big(v_{*}^{il}+ \bar{g}_{jl}\big)(t,x),\; \varpi(t,x)  \} \right\}, } \\
\end{array}
$$
where,
$$
\varpi(t,x) = - p -b(t,x)^\top q
-\frac{1}{2}\textrm{Tr}(\sigma\sigma^\top(t,x)(M -2\gamma))
-f^{ij}(t,x, [(v_{*}^{kl}(t,x))_{(k,l) \neq (i,j)},u_{\delta,
\gamma}(t,x)], \sigma^\top (t,x)q).
$$
First we note that from (\ref{upst0x0}), $\Upsilon (t_0,x_0)<0$, since
$u_{\delta, \gamma}(t_0,x_0)=v_{*}^{ij}(t_0,x_0)+\delta$.
On the other hand by the continuity of $u_{\d,\g}$, Assumptions (H1)-(H2) on
$f^{ij}$ and finally the lower semi-continuity of $v^{kl}_*$,
$(k,l)\in \gam$, we can easily check that the function $\Upsilon$ is
$usc$. Thus, for any $\eps >0$, there exists $\eta >0$ such that for
any $(t,x) \in B_\eta$ we have $\Upsilon (t_0,x_0)\geq \Upsilon
(t,x)-\epsilon$. Since $\Upsilon (t_0,x_0)<0$, then choosing $\eps$
small enough we deduce that $\Upsilon (t,x)\leq 0$ for any $(t,x)\in
B_\eta$ with an appropriate $\eta$. It follows that the function
$u_{\d,\g}$ is a viscosity subsolution on $B_\eta$ of the following system.
$$
\begin{array}{l}\min \left\{
\varrho(t,x) - \max_{k\neq i} \big(v_{*}^{kj} - \underline{g}_{ik}\big)(t,x), \;
  \max \{\varrho(t,x) - \min_{l \neq j} \big(v_{*}^{il}+ \bar{g}_{jl}\big)(t,x),\;  \right.\\ \left.\qq -
  \partial_t \varrho(t,x) -b(t,x)^\top D_x \varrho(t,x)
-\frac{1}{2}\textrm{Tr}(\sigma\sigma^\top(t,x)D^2_{xx}\varrho
(t,x))\right.\\\left.\qq\qq -f^{ij}(t,x, [(v_{*}^{kl}(t,x))_{(k,l)
\neq (i,j)},\varrho(t,x)], \sigma^\top (t,x)D_x \varrho(t,x)) \}
\right\}=0.
\end{array}
$$
Since, for any $(k,l)\in \gam$, $v^{kl}_*\le v^{kl,*}$ and since
$f^{ij}$ satisfies the monotonicity condition (H2),  $u_{\d,\g}$
is also a viscosity subsolution on $B_\eta$ of the system
\be
\label{viscointermed}
\begin{array}{l}\min \left\{
\varrho(t,x) - \max_{k\neq i} \big(v^{kj,*} -
\underline{g}_{ik}\big)(t,x); \;
  \max \{\varrho(t,x) - \min_{l \neq j} \big(v^{il,*}+ \bar{g}_{jl}\big)(t,x),\;  \right.\\\left.\qq -
  \partial_t \varrho(t,x) -b(t,x)^\top D_x \varrho(t,x)
-\frac{1}{2}\textrm{Tr}(\sigma\sigma^\top(t,x)D^2_{xx}\varrho
(t,x))\right.\\\left.\qq\qq -f^{ij}(t,x, [(v^{kl,*}(t,x))_{(k,l)
\neq (i,j)},\varrho(t,x)], \sigma^\top (t,x)D_x \varrho(t,x)) \}
\right\}=0
\end{array}
\ee
Next, as $(p,q,M)\in\mathcal{{J}}^{-}(v^{ij}_{*}(t_0, x_0))$ then
$$
\begin{array}{l} v^{ij}(t,x)\ge
 v^{ij}_*(t,x)\ge v^{ij}_{*}(t_0, x_0) +p(t-t_0)+ \langle{q, x-x_0\rangle} \\\qq\qq \qq \qq +\frac{1}{2}
 \langle{M(x-x_0), (x-x_0)
\rangle}+ o(|t-t_{0}|)+o(|x-x_{0}|^{2}).\\
\end{array} $$
In view of the definition of $u_{\delta, \gamma}$ and taking $
\delta=\frac{r^2}{8} \gamma$, it is easily seen that
$$
v^{ij}(t,x) > u_{\delta, \gamma}(t,x),
$$
as soon as $\frac{r}{2} < |x-x_0| \le r$ and $r$ small enough. We
now take $r\leq \eta$ and consider the function $\tilde
u^{ij}$:
$$\tilde u^{ij}(t,x) =\left\{  \begin{array}{ll}
 \max(v^{ij}(t,x),\; u_{\delta, \gamma}(t,x)),& \;\textrm{if} \; (t,x) \in B_{r},\\
v^{ij}(t,x), & \;\textrm{otherwise.}\\
   \end{array}\right.
$$
Then taking into account of (\ref{viscointermed}) and using Lemma 4.2 in Crandall {\it et al.} (1992) (\cite{Crandalllions92}), it follows that
$\tilde u^{ij}$ is a subsolution of (\ref{equationseparee}). Next,
as $\tilde u^{ij}\geq v^{ij}$ and using once more the monotonicity
assumption (H2) on $f^{kl}$, $(k,l)\in \gam$, we get that
$[(v^{kl})_{(k,l)\neq (i,j)}, \tilde u^{ij}]$ is also a
subsolution of (\ref{main-systvi}) which belongs to $\Pi_g$. Thus,
thanks to the Comparison Theorem \ref{comparison},
$[(v^{kl})_{(k,l)\neq (i,j)}, \tilde u^{ij}]$ belongs also to
$\cU_{m_0}$. Finally, in view of the definition of $v^{ij}_*$, there
exists a sequence $(t_n, x_n, v^{ij}(t_n, x_n))_{n\geq 1}$ that
converges to $(t_0, x_0, v_{*}^{ij}(t_0, x_0))$. This implies that
$$
\displaystyle{\lim_{n}(\tilde u^{ij} -v^{ij})(t_n, x_n) =
 \big(u_{\delta, \gamma} - v_{*}^{ij})(t_0, x_0) =v_{*}^{ij}(t_0, x_0)+ \delta -v_{*}^{ij}(t_0, x_0)
 >0.}
 $$
 This means that there are points $(t_n,x_n)$ such that $\tilde
 u^{ij}(t_n,x_n)>v^{ij}(t_n,x_n)$. But this contradicts the
 definition of $v^{ij}$, since $[(v^{kl})_{(k,l)\neq (i,j)}, \tilde
u^{ij}]$  belongs also to $\cU_{m_0}$. Therefore, $(v^{ij})_{\ij}$
is also a supersolution for (\ref{main-systvi}).

\ms\no Now, in view of Corollary \ref{unicite},
$(^{m_0}v^{ij})_{\ij}$ is the unique continuous viscosity solution in the class $\Pi_g$
of (\ref{main-systvi}). Thus, using once more uniqueness,
we deduce that $(^{m_0}v^{ij})_{\ij}$ does not depend on $m_0$.
\qed \bs

As above, let us denote by $(v^{i,j})_{\ij}$ the family
$(^{m_0}v^{i,j})_{\ij}$. Here is the second main result of the paper.

\begin{thm}\lb{mainthm2}For any $\ij$, $\bar v^{ij}=v^{ij}$, i.e.,
$(v^{ij})_{\ij}$ is continuous and is the unique viscosity solution in the class $\Pi_g$
of the system (\ref{main-systvi}).
\end{thm}
\pr. For any $\ij$ and  $m_0\in \mathbb{N}$ we have
$$
\bar v^{ij}\le v^{ij}\le \bar v^{ij,m_0}.
$$
Taking the limit as $m_0\rw \infty$ we obtain $\bar v^{ij}=v^{ij}$, for all $\ij$. Finally, Theorem \ref{main1} completes the proof. \qed

\ms\no As a by-product of this result, we have the following theorem for the family
$(\underline v^{ij})_{\ij}$.
\begin{thm}\label{mainthm3}
The functions $(\underline v^{ij})_{\ij}$ are continuous, of
polynomial growth and unique viscosity solution in the class $\Pi_g$
of the following system of variational inequalities. For every $\ij$,
\be  \lb{approxdessous}\left\{
\begin{array}{l}
\max\left\{ \underline v^{ij}(t,x)-\min_{l\in \gamj}(\underline
v^{il}+\bar g_{jl})(t,x),\right.\min \left\{ \underline
v^{ij}(t,x)- \max_{k\in \gami}(\underline v^{kj}-\underline
g_{ik})(t,x), \right. \\\left.\left.\qquad \qq -\partial_t
\underline v^{ij}(t,x)- {\cal L}\underline
v^{ij}(t,x)-f^{ij}(t,x,(\underline v^{kl}(t,x))_{(k,l)\in
\Gamma^{1} \times \Gamma^{2}}, \sigma^{\top}(t,x)D_x
\underline v^{ij}(t,x))\right\}\right\} =0\\
\underline v^{ij}(T,x)=h^{ij}(x).
\end{array}\right.\ee
\end {thm}
\pr. It is enough to consider $(-\underline v^{ij})_{\ij}$
which, in view of Theorem \ref{mainthm2}, is continuous, of polynomial growth and
the unique viscosity solution of the following
system. For all $\ij$,
\be \left\{
\begin{array}{l}
\min\left\{ v^{ij}(t,x)-\max_{l\in \gamj}(v^{il}-\bar
g_{jl})(t,x),\right.\max \left\{ v^{ij}(t,x)- \min_{k\in
\gami}(v^{kj}+\underline g_{ik})(t,x), \right.
\\\left.\left.\qquad \qq -\partial_t v^{ij}(t,x)- {\cal
L}v^{ij}(t,x)+f^{ij}(t,x,( -v^{kl}(t,x))_{(k,l)\in \Gamma^{1}
\times \Gamma^{2}}, -\sigma^{\top}(t,x)D_x
 v^{ij}(t,x))\right\}\right\} =0\\
v^{ij}(T,x)=-h^{ij}(x).
\end{array}\right.
\ee Using now the result by Barles (\cite{Barles}, pp. 18), we obtain that
$(\underline v^{ij})_{\ij}$ are continuous, of polynomial growth
and unique viscosity solution in the class $\Pi_g$ of system
(\ref{approxdessous}), which is the desired result. \qed
\begin{remark}
We do not know whether or not we have $\underline v^{ij}=\bar
v^{ij},\,\,\ij$.
\end{remark}

\section{Appendix}
\textbf{A1}. {\it Comparison of solutions of multi-dimensional BSDEs} (\cite{Hupeng}, Theorem 1, pp.135)

\begin{thm}\label{a1}
Let $(Y,Z)$ (resp. $(\bar Y,\bar Z)$) be the solution of the
$k$-dimensional BSDE associated with
$(f:=(f_i(t,\omega,y,z))_{i=1,k},\xi)$ (resp. $(\bar f:=(\bar
f_i(t,\omega,y,z))_{i=1,k}, \bar \xi)$) where:

(i) $\xi$ and $\bar \xi$ are square integrable ${\cal F}_T$-random
variables of $\R^k$ ;

(ii) the functions $f(t,\omega,y,z)$ and $\bar f(t,\omega,y,z)$
defined on $[0,T]\times \Omega \times \R^{k+k\times d}$ are
$\R^k$-valued, Lipschitz in $(y,z)$ uniformly in $(t,\omega)$ and
the process $(f(t,\omega,0,0))_{t\le T}$ (resp. $(\bar
f(t,\omega,0,0))_{t\le T}$) belongs to $\cH^{2,k}$;

(iii) for any $i=1,\ldots,k$, the $i-th$ component $f_i$ (resp. $\bar
f_i$) of $f$ (resp. $\bar f$) depends only on the $i-th$ row of the
matrix $z$.

\noindent If there exists a constant $C\geq 0$, such that for any
$y, \bar y \in \R^k$, $z, \bar z \in \R^{k\times d}$
$$-4\sum_{i=1,k}y_i^-(f_i(t,\omega,y_1^++ \bar y_1,\ldots,y_k^++ \bar y_k,z) -\bar
f_i(t,\omega,\bar y_1,\ldots,\bar y_k,\bar z))\leq 2
\sum_{i=1,k}\ind_{[y_i<0]}|z-\bar z|^2+C \sum_{i=1,k}(y_i^-)^2
$$
where $y_i^+=max(y_i,0)$ and $y_i^-=max(-y_i,0)$. Then for any
$i=1,\ldots,k$, $\P$-a.s., $Y^i\leq \bar Y^i$. \qed
\end{thm}
\bs

\noindent \textbf{A2}. {\it Systems of reflected BSDEs with one
inter-connected barrier and their related systems of variational
inequalities} (see e.g. \cite{Hamzhang} or \cite{Hammorlais11}).

Let $\cJ:=\{1,\dots,m\}$ and let us consider the following
functions: for $i,j\in \cJ$,
$$
\begin{array}{l}
f_i: \,  (t,x,y^1,\ldots,y^m,z)\in [0,T]\times \R^{k+m+d}\mapsto
f_i(t,x,y^1,\ldots,y^m,z)\,\,\in \R\,;\\
g_{ij}:\, (t,x)\in [0,T]\times\R^k\mapsto g_{ij}(t,x)\in \R\,(i\neq j);\\
h_i: x\in \R^k\mapsto  h_i(x)\in \R.\\
\end{array}$$
 We now make the following assumptions.

\ms \no {\bf [Af]}. For $i\in {\cal J}$, $f_i$ satisfies

(i) The function $(t,x)\mapsto f_i(t,x,y^1,\ldots,y^m,z)$ is continuous uniformly
w.r.t. $(\overrightarrow{y},z):=(y^1,\ldots,y^m,z)$.

(ii) The function $f_i$ is uniformly Lipschitz continuous with respect to
$(\overrightarrow{y},z):=(y^1,\ldots,y^m,z)$, i.e., for some $C\geq 0$,
$$|f_i(t,x,y^1,\ldots,y^m,z)-f_i(t,x,\bar y^1,\ldots,\bar y^m,\bar z)|\leq C(|y^1
-\bar y^1|+\dots+|y^m-\bar y^m|+|z-\bar z|).$$

(iii) The mapping $(t,x)\mapsto f_i(t,x,0,\dots,0)$ is ${\cal
B}([0,T]\times \R^k)$-measurable and of polynomial growth i.e. it
belongs to $\Pi_g$;\ms

(iv) \underline{\it{Monotonicity}}. $\forall i\in {\cal J}$, for any
$k\in {\cal J}^{-i}$, the mapping $y_k\in \R\mapsto
f_i(t,x,y_1,\ldots,y_{k-1},y_k,y_{k+1},\ldots,y_m)$  is non-decreasing
whenever the other components
$(t,x,y_1,\ldots,y_{k-1},y_{k+1},\ldots,y_m)$ are fixed. \bs

\noindent {\bf [Ag]}. (i) The function $g_{ij}$ is jointly continuous in $(t,x)$,
non-negative, i.e., $g_{ij}(t,x)\geq 0$, $\forall (t,x)\in
[0,T]\times \R^k$ and belongs to $\Pi_g$.

(ii) {\it The no free loop property}. for any $(t,x)\in [0,T]\times
\R^k$ and for any sequence of indexes $i_1,\ldots,i_k $ such that
$i_1=i_k$ and $card\{i_1,\ldots,i_k\}=k-1$ we have
$$g_{i_1 i_2}(t,x)+g_{i_2 i_3}(t,x)+\dots+g_{i_{k-1} i_k}(t,x)+g_{i_k i_1}(t,x)>0,
\,\,\forall (t,x)\in [0,T]\times \R^k.$$
As a convention we assume hereafter that $g_{ii}(t,x)=0$ for any $(t,x)\in \esp$ and $\ij$.\ms

\no {\bf [Ah1]}. $h_i$ is continuous, belongs to $\Pi_g$ and
satisfies:
$$\forall x\in\R,\,\,
h_i(x)\geq \max_{j\in {\cal J}^{-i}}(h_j(x) -g_{ij}((T,x)).$$ \ms

\no {\bf [Ah2]}. The function $h_i$ is continuous, belongs to $\Pi_g$ and
satisfies
$$\forall x\in\R,\,\,
h_i(x)\geq \min_{j\in {\cal J}^{-i}}(h_j(x)+g_{ij}((T,x)).
$$
Then, we have
\begin{thm}\label{a2} Assume that [Ah], [Ag] and [Ah1] are fulfilled. Then, there exist $m$ triples of
processes\\
$((Y^{i;t,x},Z^{i;t,x},K^{i;t,x}))_{i\in {\cal J}}$ that satisfy:
$\forall \,\,i\in {\cal J}$,
\be\label{systemapend}\left\{\begin{array}{l} Y^i, K^i \in \cS^2,
\,\,Z^i \in \cH^{2,d},\,\,K^i \mbox{
non-decreasing and }K^i_0=0 ;\\
Y^i_s=h_i(X^{t,x}_T)+\int_s^Tf_i(r,X^{t,x}_r,
Y^1_r,\dots,Y^m_r,Z^i_r)dr+K^i_T-K^i_s-\int_s^TZ^i_rdB_r,\,\,\forall\,\,s\leq
T\\
Y^i_s\geq \max_{j\in {\cal J}^{-i}}\{Y^j_s-g_{ij}(s,X^{t,x}_s)\},
\,\,\forall s\leq T\\
\int_0^T(Y^i_s-\max_{j\in {\cal
J}^{-i}}\{Y^j_s-g_{ij}(s,X^{t,x}_s)\})dK^i_s=0.\end{array}\right.\ee

Moreover there exist $m$ deterministic functions $(v^i(t,x))_{i\in
\cJ}$ continuous and belonging to $\Pi_g$ such that:
$$ \forall s\in [t,T], Y^{i;t,x}_s=v^i(s,X^{t,x}_s).$$ Finally
$(v^i(t,x))_{i\in \cJ}$  is the unique solution, in the sub-class of
$\Pi_g$ of continuous functions, of the following system of
variational inequalities with inter-connected obstacles: $\forall
\,\,i\in {\cal J}$
\be \label{sysvi01apend} \left\{
\begin{array}{l}
\min\left \{v_i(t,x)- \max\limits_{j\in{\cal
J}^{-i}}(-g_{ij}(t,x)+v_j(t,x)),\right.\\\left.
\qquad\qquad-\partial_tv_i(t,x)- {\cal
L}v_i(t,x)-f_i(t,x, v^1(t,x), \dots,v^m(t,x), \sigma^\top(t,x)D_xv^i(t,x))\right\}=0\,\,;\\
v_i(T,x)=h_i(x).
\end{array}\right.
\ee
\end{thm}
\begin{remark}\label{rmkappendice}
In equations (\ref{systemapend}) and (\ref{sysvi01apend}), if instead
we have required an upper barrier reflection, then one would have
obtained a similar result which can be stated as follows. \bs

\noindent Assume that [Af], [Ag] and [Ah2] are fulfilled. Then there
exist $m$ triples of processes $((\tilde Y^{i;t,x},\tilde
Z^{i;t,x},\tilde K^{i;t,x}))_{i\in {\cal J}}$ that satisfy, for all $i\in {\cal J}$,
\be\label{systemapend2}\left\{\begin{array}{l}
\tilde Y^i, \tilde K^i \in \cS^2, \,\,\tilde Z^i \in
\cH^{2,d},\,\,\tilde K^i \mbox{
non-decreasing and }\tilde K^i_0=0 ;\\
\tilde Y^i_s=h_i(X^{t,x}_T)+\int_s^Tf_i(r,X^{t,x}_r, \tilde
Y^1_r,\dots,\tilde Y^m_r,\tilde Z^i_r)dr-(\tilde K^i_T-\tilde
K^i_s)-\int_s^T\tilde Z^i_rdB_r,\,\,\quad t\le s\leq
T,\\
\tilde Y^i_s\leq \min_{j\in {\cal J}^{-i}}\{\tilde
Y^j_s+g_{ij}(s,X^{t,x}_s)\},
\,\,\quad t\le s\leq T, \\
\int_0^T(\tilde Y^i_s-\min_{j\in {\cal J}^{-i}}\{\tilde
Y^j_s+g_{ij}(s,X^{t,x}_s)\})d\tilde K^i_s=0.\end{array}\right.
\ee
Moreover, there exist $m$ deterministic functions $(\tilde
v^i(t,x))_{i\in \cJ}$ continuous and belong to $\Pi_g$ such that:
$$
\tilde Y^{i;t,x}_s=\tilde v^i(s,X^{t,x}_s),\qquad s\in [t,T].
$$
Finally, $(\tilde v^i(t,x))_{i\in \cJ}$  is the unique solution, in the
subclass of $\Pi_g$ of continuous functions, of the following
system of variational inequalities with interconnected obstacles. For all $i\in {\cal J}$
\be \label{sysvi01apend2} \left\{
\begin{array}{l}
\max\left \{\tilde v_i(t,x)- \min\limits_{j\in{\cal
J}^{-i}}(g_{ij}(t,x)+\tilde v_j(t,x)),\right.\\\left.
\qquad\qquad-\partial_t\tilde v_i(t,x)- {\cal
L}\tilde v_i(t,x)-f_i(t,x, \tilde v^1(t,x), \dots,\tilde v^m(t,x), \sigma^\top(t,x)D_x\tilde v^i(t,x))\right\}=0\,\,;\\
\tilde v_i(T,x)=h_i(x).
\end{array}\right.
\ee
The proof of this result is obtained straightforward from
Theorem \ref{a2}, in considering the equations satisfied by
$((-\tilde Y^i,-\tilde Z^i, \tilde K^i))_{i\in \cJ}$.
\end{remark}

\ms\no\textbf{A3}. {\it Linearization procedure of Lipschitz
functions.}
Let $f$ be a function from $\R^2$ to $\R$ which with $(x_1,x_2)$
associates $f(x_1,x_2)$ which is Lipschitz in its arguments. Then, we
can write \be\lb{linproc}
\begin{array}{ll}
f(x_1,x_2)-f(y_1,y_2)&=f(x_1,x_2)-f(y_1,x_2)+f(y_1,x_2)-f(y_1,y_2)\\
{}&=\ind_{x_1-y_1\neq
0}\frac{f(x_1,x_2)-f(y_1,x_2)}{x_1-y_1}(x_1-y_1)+ \ind_{x_2-y_2\neq
0}\frac{f(y_1,x_2)-f(y_1,y_2)}{x_2-y_2}(x_2-y_2)\\
{}&=a_1(x_1,x_2,y_1).(x_1-y_1)+ a_2(x_2,y_1,y_2).(x_2-y_2)
\end{array}
\ee where, $a_1$ and $a_2$ are measurable functions and bounded,
i.e.,
$$
|a_1(x_1,x_2,y_1)|\vee |a_2(x_2,y_1,y_2)|\le C(f), \,\,i=1,2,
$$
where, $C(f)$ is the Lipschitz constant of $f$. Moreover, if $f$ is
non-decreasing with respect to $x_1$ (resp. $x_1$) when $x_2$ (resp.
$x_1$) is fixed, then $a_i\geq 0$, $i=1,2$.

\ms\no Linearizing $f$ consists of writing $f$ as
$$
f(x_1,x_2)-f(y_1,y_2)=a_1(x_1,x_2,y_1).(x_1-y_1)+
a_2(x_2,y_1,y_2).(x_2-y_2).
$$
\noindent \textbf{A4}. {\it Representation of a penalization scheme of
two barriers reflected BSDE.}

\ms\no For $n\geq 0$ let $(Y^n, Z^n,K^{+,n})$ be the solution of the
following one barrier reflected BSDE.
 \be \label{schemapenalisation}
\left\{
\begin{array}{l}
(Y^n,Z^n,K^{+,n})\in \cS^{2,1}\times \cH^{2,d}\times \cA^{+,2}\\
Y^n_t=\xi+\int_t^Tg(s)ds-n\int_t^T(Y^n_s-U_s)^+ds+K^{+,n}_T-K^{n,+}_t-\int_t^TZ^n_sdB_s,\quad 
t\leq T; \\Y^n_t\geq L_t,\quad 0\le t\leq T\,\,\mbox{ and
}\,\,\int_0^T(Y^n_s-L_s)dK^{+,n}_s=0 \ea \right.
\ee
where, the processes $L$ and $U$ belong to $\cS^{2,1}$, $(g(s))_{s\leq T} \in \cH^{2,1}$,
$\xi$ is square integrable and $\cF_T$-measurable. Moreover, we
require that $L\leq U$ and $L_T\le \xi$. Under these conditions, the
solution $(Y^n, Z^n,K^{+,n})$ exists and is unique (see e.g.
\cite{elkarouietal}). Next, for $n\ge 0$ and $t\leq T$, set
$K^{n,-}_t=n\int_0^t(Y^n_s-U_s)^+ds$. Then, $K^{-,n}\in \cA^{+,2}$
and\\ $\int_0^T(Y^n_s- Y^n_s\vee U_s)dK^{-,n}_s=0$. Therefore,
the equation (\ref{schemapenalisation}) can be expressed as a BSDE with two
reflecting barriers in the following manner. For all $t\leq T$,
\be
\label{schemapenalisation2BSDE} \left\{
\begin{array}{l}
Y^n_t=\xi+\int_t^Tg(s)ds+(K^{+,n}_T-K^{+,n}_t)-
(K^{-,n}_T-K^{-,n}_t)-\int_t^TZ^n_sdB_s; \\L_t\leq Y^n_t\leq
Y^n_t\vee U_t, \\\int_0^T(Y^n_s-L_s)dK^{+,n}_s=\int_0^T(Y^n_s-
Y^n_s\vee U_s)dK^{-,n}_s=0. \ea \right.
\ee
Thus, a result by Cvitanic and Karatzas \cite{CK96} or Hamad\`ene and Lepeltier \cite{HL95} allows to
represent $Y^n$ as a value function of a Dynkin game, i.e., it holds
true that for any $t\leq T$,
$$\ba{ll}
Y^n_t&=\esssup_{\sigma \ge t}\essinf_{\t \geq
t}\E[\int_t^{\si\wedge \t}g(s)ds +L_\si \ind_{[\si < \t
]}+(Y^n_\t\vee U_\t)\ind_{[\t \leq \si<T]}+\xi
\ind_{[\t=\si=T]}|\cF_t]\\{}&= \essinf_{\t \geq
t}\esssup_{\sigma \ge t}\E[\int_t^{\si\wedge \t}g(s)ds +L_\si
\ind_{[\si < \t ]}+(Y^n_\t\vee U_\t)\ind_{[\t \leq \si<T]}+\xi
\ind_{[\t=\si=T]}|\cF_t].\ea
$$
where $\t$ and $\si$ are $\cal F$-stopping times.
\end{document}